\documentclass[journal]{IEEEtran}

\usepackage{amsmath}
\usepackage{amsthm}
\usepackage{amssymb}

\usepackage{algorithmic}
\usepackage{algorithm}

\usepackage{amsfonts,mathrsfs}

\usepackage{mathptmx}       
\usepackage{helvet}         
\usepackage{courier}        
\usepackage{type1cm}        
\usepackage{makeidx}         
\usepackage{graphicx}        
\usepackage{color}

\makeindex             

\newcommand{\AHLS}{LHS\ }

\newcommand{\AHLSS}{LHSs\ }
\newcommand{\SAHLS}{LHSs\ }

\newcommand{\Time}{\mathbb{R}_{+}}
\newcommand{\GammaS}{(\Gamma \times \Time)}
\newcommand{\GammaTSN}{\Gamma_{\mathrm{timed}}^{\infty}}
\newcommand{\GammaTS}{\Gamma_{\mathrm{timed}}^{\infty}}

\newcommand{\cQ}{\mathscr{Q}}
\newcommand{\cP}{\mathscr{P}}
\newcommand{\bU}{\mathbf{U}}
\newcommand{\bS}{\mathbf{S}}
\newcommand{\bV}{\mathbf{V}}
\newcommand{\bLa}{\Lambda}

\newcommand{\bfz}{{\mathbf 0}}

\newtheorem{Definition}{Definition}

\newtheorem{Theorem}{Theorem}

\newtheorem{Lemma}{Lemma}

\newtheorem{Remark}{Remark}

\newtheorem{Notation}{Notation}

\newtheorem{Procedure}{Procedure}

\newcommand{\HYBINP}{\mathbf{U}}

\newcommand{\bx}{x}
\newcommand{\by}{y}
\newcommand{\bu}{u}
\newcommand{\bq}{q}

\begin{document}
%
\title{Model reduction of linear hybrid systems}
%
%
%



\author{Ion Victor~Gosea,~Mihaly~Petreczky,~John~Leth,~Rafael~Wisniewski
        and~Athanasios~C.~Antoulas
\thanks{I. V. Gosea is with Data-Driven System Reduction and Identification (DRI) Group, Max Planck Institute for Dynamics of Complex Technical Systems, Sandtorstrasse 1, 39106, Magdeburg, Germany, e-mail: gosea@mpi-magdeburg.mpg.de.}
\thanks{M. Petreczky is with Centre de Recherche en Informatique, Signal et Automatique de Lille (CRIStAL),  UMR CNRS 9189, CNRS, Ecole Centrale de Lille, France, e-mail: mihaly.petreczky@ec-lille.fr}
\thanks{J. Leth and R. Wisniewski are with Department of Electronic Systems, Automation and Control, Aalborg University, Fredrik Bajers Vej 7C, Aalborg, Denmark, e-mail: jjl@es.aau.dk, raf@es.aau.dk.}
\thanks{A.C. Antoulas is with Department of Electrical and Computer Engineering,	Rice University, 6100 Main St, MS-366, Houston, TX 77005, USA, and DRI Group, Max Planck Institute for Dynamics of Complex Technical Systems, Sandtorstrasse 1, 39106, Magdeburg, and Baylor College of Medicine, 1 Baylor Plaza, Houston, TX 77030, e-mail: aca@rice.edu}}
\maketitle

\begin{abstract}
The paper proposes a model reduction algorithm for linear hybrid systems, i.e., hybrid systems with externally induced discrete events, with linear continuous subsystems, and 
linear reset maps. The model reduction algorithm is based on balanced truncation. Moreover, the paper also proves an analytical error bound for the 
difference between the input-output behaviors of the original and the reduced order model. This error bound is formulated in terms of singular values of the Gramians used for model reduction. 
\end{abstract}


%

\section{Introduction}

In this paper we propose a model reduction method for linear hybrid systems with external switching. A linear hybrid system is a hybrid system continuous states of which are governed by linear differential equations, the reset maps are linear, and the discrete-events are external inputs. 
Linear hybrid systems can be viewed as a generalization of linear switched systems \cite{D:Lib,BookSun11}, but in contrast to linear switched systems we allow
state jumps and the change of discrete states is supposed to follow the transition structure of a Moore automaton. Linear hybrid systems occur in several applications, and a well known class of piecewise-affine systems is directly related to linear hybrid systems, as the former can be viewed as a feedback interconnection of the latter with a discrete-event generator. 
The model reduction method we propose is based on balanced truncation, performed for each linear subsystem. The corresponding Gramians have to satisfy certain linear matrix inequalities (LMIs).  In addition to the novel algorithm, we propose an analytic error bound for the difference between the input-output behaviors of the original and the reduced-order models. This error bound is a direct counterpart of the well-known error bound for balanced truncation of linear systems \cite{ModRed1}, and it involves the singular values of the Gramians.

To the best of our knowledge, the contribution of the paper is new. 
Indeed, the existing methods for model reduction of hybrid
systems can be grouped into the following categories. 

\textbf{LMI-based methods}
These methods compute the matrices of the reduced order model
by solving a set of LMIs. 
           The disadvantage is that the proposed conditions are only sufficient, and the trade-off between the dimension of the reduced model and the error bound is not clear. Moreover, the computational complexity of solving those LMIs might be too high. 
  Without claiming completeness, we mention the following papers 
 \cite{glw06,zbs09,zsbw08,zcw14}. 
First of all, the cited papers do not deal with linear reset maps. 
Moreover, in contrast to the cited papers, the current paper proposes a method, whose applicability depends on the
           existence of solution for a few simple LMIs which are necessary to find the observability/controllability Gramians.
            Once the existence of these Gramians is assured, the model reduction method can be applied. Moreover, there is an analytic error 
            bound and the trade-off between 
            the approximation error and the dimension of the reduced system is formalized in terms of the singular values of those Gramians.         

\textbf{Methods based on local Gramians}
  

The algorithms which belong to this class are based on finding observability/controllability Gramians for each
linear subsystem. 
They are solutions of LMIs derived by relaxing the classical Lyapunov-like equations for observability/controllability Gramians. 
The disadvantage of these methods is that often there are no error bounds or the reduced order model need not be well-posed.  
Examples of such papers include \cite{mtc12,pp16,bgmb10,bmb11,bmb12,gpaf18}.
Note that to the best of our knowledge, the only algorithm which always yields a well-posed
linear switched system of the same type as the original one and for which there exists an analytic error bound is the one of
\cite{gpaf18}. Even this algorithm provides an error bound only for sufficiently slow switching signals (i.e., switching sequences with a suitable minimal dwell time).
The method of this paper is an extension of \cite{gpaf18}. The main difference between the current paper and \cite{gpaf18} is the following: 
\begin{itemize}
\item
      In contrast to  \cite{gpaf18}, the error bound of this paper
      no longer uses the assumption of minimum dwell time. However, this comes at price, as the LMIs involved are more conservative. 

\item The discrete states are no longer assumed to be inputs, but
      they are states of the system and they are assumed to evolve according to a 
      Moore-automaton.  However, the Moore-automaton is driven by
      discrete events which are external inputs. That is, the system class
     considered in this paper is more general than that of \cite{gpaf18}.
\end{itemize}

More recently, a balancing truncation method for linear switched systems that are characterized by constrained
switching scenarios was proposed in \cite{gpba20}. The technique is based on defining generalized Gramians for each discrete mode, specifically tailored to particular switching scenarios.

%
%
%

\textbf{Methods based on common Gramians}
   These methods rely on finding the same observability/controllability Gramian for each linear subsystem. In most contributions, the Gramians
   are derived as solutions of a suitable LMI.  
   Such algorithms were described in \cite{sw11,sw12} and an analytic error bound was derived in \cite{PetreczkyNAHS}.
   The results of this paper can also be viewed as a direct extension of
   \cite{PetreczkyNAHS}. In particular, when applied to a linear switched system
   of the type studied in \cite{PetreczkyNAHS}, the results of the present paper
   boil down to those of \cite{PetreczkyNAHS}. With respect to \cite{PetreczkyNAHS}, the main novelty of the present paper is that it considers a system class which is much larger than the one of \cite{PetreczkyNAHS}. 
   Nevertheless, some methods that do not rely on solving LMIs are also available. For example, in \cite{su18} a balancing procedure based on  recasting the original linear switched system as an envelope linear time-invariant system with no switching was proposed. Additionally, a balancing procedure based on reformulating the original system as a bilinear system with no switching was presented in \cite{pgb18}.

\textbf{Moment matching}
       The idea behind these algorithms is to find a reduced order switched system such that certain coefficients of the series expansions of
       the input-output maps of the original and the reduced order system coincide. The series expansion can be the Taylor series with respect to switching times, in
       which case the so-called Markov parameters are matched. Alternatively, the series expansion can be a Laurent-series expansion of a multivariate Laplace transform
       of the input-output map around a certain frequency. The former approach was pursued in \cite{MertBastug:Automatica,MertBastug:TAC,MertBastugPhd} , the   latter in \cite{VictorLoewner}. While 
            those methods do not allow for analytical error bounds, under suitable assumption it can be guaranteed that the reduced
            model will have the same input-output behavior for certain switching signals \cite{MertBastug:Automatica,MertBastug:TAC,MertBastugPhd}.
     A somewhat different approach is that of \cite{Astolfi}, which 
          considers  switched systems with autonomous switching and it proposed a model reduction procedure which 
            guarantees that the reduced model has the same steady-state output response to certain inputs as the original model. 

The results of the present paper are based on 
balanced truncation. As a result, in contrast to the cited papers, we are able to propose an analytic error bound. Moreover, the class of systems considered in this paper is much larger than that of the cited papers. In particular, we allow reset maps and the evolution of the discrete states is governed by a Moore-automaton.

The paper is structured as follows. In Section \ref{petreczky:sect2} we  fix the notation and we present the formal definition of linear hybrid systems and of some related concepts. In Section \ref{bal:trunc} we present a balanced truncation algorithm for model reduction and an analytical error bound for this algorithm.
In Section \ref{sect:num} we present a numerical example to illustrate the proposed algorithm. 
In Appendix \ref{tech:proofs} we present the proofs of the technical results used in the paper.

\section{Preliminaries}
\subsection{Notation}\label{sect2}

Let $\mathbb{N}$ denote the set of natural numbers including $0$, and $\mathbb{R}_{+}=[0, +\infty)$ denote the positive \emph{real time-axis}. We denote by $PC(A,B)$ the set of all \emph{piecewise-continuous maps} $A\to B$, and by $L_{2}(A,B)$ the set of all \emph{Lebesgue measurable maps} $A
\to B$. The $L_2$-norm and Euclidean 2-norm are denoted by $\|\cdot\|_{L_2}$ and $\|\cdot\|_{2}$ respectively.

\subsection{ Linear hybrid systems: definition and basic concepts }
\label{petreczky:sect2}
\begin{Definition}[\AHLS]
	A \emph{linear hybrid system}\index{linear hybrid system} $H$ (abbreviated as \AHLS) \index{\AHLS}
	is a tuple
	\begin{align} \label{petreczky:ahlsd}
	\begin{split}
	H = &	(Q,\Gamma,O,\delta,\lambda, \{n_q,A_q,B_q,C_q \}_{q \in Q}, \\
	 & \{ M_{q_{1},\gamma,q_{2}} \}_{q_{2} \in Q, \gamma \in \Gamma, q_1=\delta(q_2,\gamma)}, h_0),
	 	\end{split}
	\end{align}
	where
	\begin{enumerate}
		\item $Q$ is a finite set, called \emph{the set of discrete states},
		\item $\Gamma$ is a finite set, called \emph{the set of discrete events},
		\item $O$ is a finite set, called \emph{the set of discrete outputs},
		\item $\delta:Q \times \Gamma \rightarrow Q$ is a function called
		the \emph{discrete state-transition map},
		\item $\lambda:Q \rightarrow O$ is a function called the
		\emph{discrete readout map}.
		\item
		$\Sigma_q=(A_q,B_q,C_q)$, $q \in Q$ is 
		\emph{the linear system in the discrete state $q$} and 
		\( A_q \in \mathbb{R}^{n_q \times n_q}, B_q \in \mathbb{R}^{n_q \times m},
		C_q \in \mathbb{R}^{p \times n_q} \)
		are the matrices of this linear system.
		\item
		$M_{q_1,\gamma,q_2} \in \mathbb{R}^{n_{q_1} \times n_{q_2}}$ are matrices for 
		all  $q_2 \in Q, \gamma \in \Gamma, q_1=\delta(q_2,\gamma)$, which
		are called \emph{reset maps}.
		\item
		$h_0=(q_0,x_0)$ is the \emph{initial state}, where $q_0 \in Q$ and 
		$x_0 \in \mathbb{R}^{n_{q_0}}$.
		\end{enumerate}
		The space $\mathbb{R}^{n_q}$, $q \in Q$, $0 < n_q \in \mathbb{N}$,
		is called the \emph{continuous state space associated with the discrete state
			$q$}, $\mathbb{R}^m$ is called the
		\emph{continuous input space}, $\mathbb{R}^{p}$ is called
		the \emph{continuous output space}.
		The \emph{state space $\mathcal{H}_{H}$ of $H$} is the set
		\( \mathcal{H}_{H}=\bigcup_{q \in Q} \{q\} \times \mathbb{R}^{n_q} \). 
		\end{Definition}
		
		\begin{Notation}
		An element $x \in \mathcal{H}_{H}$  comprises of a  pair $x = (q,x_q)$ with $q \in Q$ and $x_q \in \mathbb{R}^{n_q}$. In many places in the article, we will suppress the notation and write $ x = x_q$, when it is clear from the contents which discrete mode $x$ is in.
		\end{Notation}
		Notice that the linear control systems associated
		with different discrete states may have different state-spaces, but they have the same
		input and output space.
		The intuition behind the definition of a  linear hybrid system is as follows.
		We associate a linear system 
		\begin{align}
		\label{petreczky:lti:hyb}
		\Sigma_q \
		\begin{cases}
		\dot \bx = A_{\bq}\bx + B_{\bq}\bu \\
		\by = C_{\bq}\bx
		\end{cases},
		\end{align}
		with each discrete state $q \in Q$. As long as we are in the
		discrete state $q$, the state $x$ and the continuous
		output $y$ develops according to \eqref{petreczky:lti:hyb}. 
		The discrete state can change only if a discrete event $\gamma \in \Gamma$
		takes place. If a discrete event $\gamma$ occurs at time $t$, then 
		the new discrete state $q^{+}$ is determined by applying the
		discrete state-transition map
		$\delta$ to $q$, i.e.
		\( q^{+}=\delta(q,\gamma) \).
		The new continuous-state $x^{+}(t) \in \mathbb{R}^{n_{q^{+}}}$ 
		is computed from the current continuous
		state $x(t^-)=\lim_{s\uparrow t}x(s)$ by applying the \emph{reset map}
		$M_{q^{+},\gamma,q}$ to $x(t^-)$, i.e.
		\(  \bx^{+}(t)=M_{q_{+},\gamma,q}\bx(t^-) \).
		After the transition, the continuous state $\bx$ and the continuous 
		output $\by$ evolve
		according to the linear system associated with the new discrete state $q^{+}$,
		started from the initial state $x^{+}(t)$.
		Finally, when in a discrete state $q \in Q$, the system produces
		a discrete output $o = \lambda(q)$.

		Notice that \emph{the discrete events are external inputs.
			All the continuous subsystems
			are defined with the same inputs and outputs, but on possibly different
			state-spaces.}
		Below we will formalize the intuition described above, by 
		defining input-to-state and input-output maps for \AHLS.  To this end,	we need the following. 
		
		\begin{Definition}[Timed sequences]
		A \emph{timed sequence of discrete events}\index{timed sequences of discrete event} is an  infinite sequence over
		the set $\GammaS$, i.e. it is a sequence of the form
		\begin{equation}
		\label{petreczky:def:tim:seq}
		w=(\gamma_{1},t_{1})(\gamma_{2},t_{2}) \cdots (\gamma_{k},t_{k}) \cdots ,
		\end{equation}
		where $\gamma_{i}\in \Gamma$, $k > 0$ are discrete events,  and
		$t_{i} \in  \Time$ are time instances, and $\lim_{k \rightarrow \infty} \sum_{i=1}^{k} t_i = \infty$. 
		We denote the set of timed sequences  of discrete events by 
		$\GammaTSN$. 
		\end{Definition}
		The interpretation of a timed sequence $w \in \GammaTS$ as 
		above is the following.
		If $w$ is of the form \eqref{petreczky:def:tim:seq}, then 
		$w$ represents the scenario, when
		the event $\gamma_{i}$ took place \emph{after} the event 
		$\gamma_{i-1}$ and $t_{i}$ \emph{is the time which has passed between
			the arrival of $\gamma_{i-1}$ and the arrival of $\gamma_{i}$}, i.e.
		$t_{i}$ is the difference of the arrival times of 
		$\gamma_{i}$ and $\gamma_{i-1}$. Hence, $t_{i} \ge 0$ but
		we allow $t_{i}=0$, i.e., we allow $\gamma_{i}$ to arrive
		instantly after $\gamma_{i-1}$. If $i=1$, then $t_{1}$ is simply
		the time when the first event $\gamma_{1}$ arrived.

		

		\begin{Notation}[Inputs $\HYBINP$]
		Denote by $\HYBINP=L_2(\mathbb{R}_{+},\mathbb R^m) \times \GammaTS$
		the set of inputs of a \AHLS.
		\index{$\HYBINP$}
		\end{Notation}
		
		If $(u,w) \in \HYBINP$, then $u$ represents the continuous-valued input to be fed to
		the system, $w$ represents the timed-event sequence.
		Below we define the notion of input-to-state and input-output maps for
		\AHLSS. 
		These functions map elements from $\HYBINP$ to states and outputs respectively.

		In the rest of this section, \emph{$H$ denotes a \AHLS of the form \eqref{petreczky:ahlsd}}.
		\begin{Definition}[Input-to-state map]
		\label{petreczky:AHLS:def1}
		The \emph{input-to-state map\index{input-to-state map for linear hybrid systems} of $H$ induced by the initial
			state $h_{0}=(q_{0},x_{0}) \in \mathcal{H}_{H}$ of $H$} 
		is the function
		$\xi_{H,h_{0}}:\HYBINP \rightarrow PC(\mathbb{R}_{+},\mathcal{H}_H)\times PC(\mathbb{R}_{+},Q)$
		 such that the following holds.
		For any $(u,w) \in \HYBINP$, where $w$ is of the form \eqref{petreczky:def:tim:seq}, define $T_0=0, T_i=\sum_{j=1}^{i} t_j$, $i \in \mathbb{N}$. Then 
		$\xi_{H,h_0}(u,w)=(\bx,\bq)$ such that
		\begin{enumerate}
		\item $\bq(t)=q_i$, $t \in [T_i,T_{i+1})$, where $q_0=q_I$ and $q_{i+1}=\delta(q_i,\gamma_{i+1})$  for all $i \in \mathbb{N}$
		\item 
		The restriction of $\bx$ to $[0,T_1)$
		is the unique solution (in the sense of Caratheodory) of the differential equation $\dot z(t) =A_{q_I}z(t)+B_{q_I}u(t)$, $z(0)=x_I$ on $[0,T_1)$, and
		the restriction of $\bx$ to $[T_i,T_{i+1})$ for $i > 0$ is the unique solution (in the sense of Caratheodory) of the differential equation
		$\dot z(s) =A_{q_i}z(s)+B_{q_i}u(s)$, $z(T_{i})=M_{q_{i+1},\gamma_{i+1},q_i}\lim_{t \uparrow T_{i}} \bx(t)$.
		\end{enumerate}
		\index{$\xi_{H,h_0}$}
		\end{Definition}
		\begin{Definition}[Input-output map]
		\label{petreczky:AHLS:def2}
		The \emph{input-output map  of the system $H$ \index{input-output map of a linear hybrid system}
			induced by the state $h \in \mathcal{H}_{H}$ of $H$} is the function
		$\upsilon_{H,h}:\HYBINP \rightarrow PC(\mathbb{R}_{+},O)\times PC(\mathbb{R}_{+},\mathbb R^p)$
		 defined
		as follows: for all $(u,w) \in \HYBINP$, 
		$\upsilon_{H,h}(u,w)=(\mathbf{o},\by)$, such that if
		$(\bq,\bx)=\xi_{H,h}(u,w)$, then
		\[ \mathbf{o}(t)=\lambda(\bq(t)), ~ \by(t)=C_{\bq(t)}\bx(t). \]
		The input-output map $\upsilon_{H,h}$ induced by
		the initial state $h_0$ is called the \emph{input-output map} of $H$
		and it is denoted by $\upsilon_{H}$.
		\index{$\upsilon_{H,h}$}\index{$\upsilon_H$}
		\end{Definition}

		\section{Balanced truncation}
\label{bal:trunc}
		Consider an \AHLS $H$ of the form \eqref{petreczky:ahlsd} with initial condition $h_0=(q_0,x_0)$ such that $x_0=0$. 
		\begin{Definition}
		\label{obs:gram:def}
		A collection $\{\cQ_q\}_{q \in Q}$ of positive definite matrices is called a collection of generalized observability Gramians of $H$, if for all $q \in Q$,
		\begin{equation}
		\label{obs:gram:eq}
		\begin{split}
		& A_q^T\cQ_q+\cQ_qA_q+C_q^TC_q < 0,  \\
		& \forall \gamma \in \Gamma,~ q^{+}=\delta(q,\gamma): 
		M_{q^+,\gamma,q}^T\cQ_{q^+} M_{q^+,\gamma,q} - \cQ_{q} \leqslant 0.
		\end{split}
		\end{equation}
		\end{Definition}
		\begin{Definition}
		\label{contr:gram:def}
		A collection $\{\cP_q\}_{q \in Q}$ of positive definite matrices is called a collection of generalized reachability Gramians of $H$, if for all $q \in Q$,
		\begin{equation}
		\label{contr:gram:eq}
		\begin{split}
		& A_q\cP_q+\cP_qA_q^T+B_qB_q^T < 0,  \\
		& \forall \gamma \in \Gamma,~ q^{+}=\delta(q,\gamma): 
	M_{q^+,\gamma,q}\cP_{q} M_{q^+,\gamma,q}^T - \cP_{q^{+}} \leqslant 0 .
		\end{split}
		\end{equation}
		\end{Definition}
		\begin{Remark}
		\label{LMI:rem1}
		The LMIs in \eqref{obs:gram:eq} can be rewritten as follows
		\begin{equation}
		\label{obs:gram:eq2}
		\begin{split}
		& \forall x \in \mathbb{R}^{n_q} : 2 (A_qx)^T\cQ_qx \le -\|C_qx\|^2_2, \\ 
		 &  x^TM_{q^+,\gamma,q}^T\cQ_{q^+}M_{q^+,\gamma,q}x \le x^T\cQ_qx.
		\end{split}
		\end{equation}
		The LMIs in \eqref{contr:gram:eq} can be rewritten as follows
		\begin{equation}
		\label{contr:gram:eq2}
		\begin{split}
		\forall x \in \mathbb{R}^{n_q},u \in \mathbb{R}^m:~ & 2(A_qx+B_qu)^T\cP_q^{-1}x \leq \|u\|^2_2, \\
		& \hspace{-24mm}  x^TM_{q^+,\gamma,q}^T\cP_{q^+}^{-1}M_{q^+,\gamma,q}x \le x^T\cP_q^{-1}x.
		\end{split}
		\end{equation}
		\end{Remark}
		\begin{Definition}
		We say that the \AHLS $H$ is quadratically stable, if there exists a collection $P_q > 0$, $q \in Q$, such that
		\begin{equation}
		\label{stab:gram:eq}
		\begin{split}
		& A_q^TP_q+P_qA_q  < 0,  \\
		& \forall \gamma \in \Gamma,~ q^{+}=\delta(q,\gamma): 
		M_{q^+,\gamma,q}^TP_{q^+} M_{q^+,\gamma,q} - P_{q} \leq 0.
		\end{split}
		\end{equation}
		\end{Definition}
	\noindent
	Next, we will briefly sketch the proof for the fact that the LMIs in \eqref{contr:gram:eq} are equivalent to those in \eqref{contr:gram:eq2}. In what follows we  use the following classical result.
			\begin{Lemma}
			Assume P and Q are negative definite matrices, i.e., $P, Q <0$. Then it follows that
			\begin{equation}
			\left[ \begin{matrix}
			P & A \\ B & Q
			\end{matrix} \right] \leq 0   \Leftrightarrow P - A Q^{-1} B \leq 0.
			\end{equation}
			\end{Lemma}
			\noindent
			Hence, using the above lemma, one can write that
                      \begin{equation}
                        \label{contr:gram:eq2:pf}
			\begin{split}
			& M_{q^+,\gamma,q}P_qM_{q^+,\gamma,q}^T - P_{q^+} \leq 0  \Leftrightarrow   \left[ \begin{matrix}
			-P_{q^+} & M_{q^+,\gamma,q} \\ M_{q^+,\gamma,q}^T & - P_{q}^{-1}
			\end{matrix} \right]   \leq 0  \\
			&  \Leftrightarrow \left[ \begin{matrix}
			-P_{q}^{-1} & M_{q^+,\gamma,q}^T \\ M_{q^+,\gamma,q} & -P_{q^+}
			\end{matrix} \right] \leq 0  \Leftrightarrow  -P_q^{-1} +M_{q^+,\gamma,q}^TP_{q^+}^{-1}M_{q^+,\gamma,q} \leq 0.
			\end{split}
                      \end{equation}
		This immediately shows that the second inequality in \eqref{contr:gram:eq2} holds for any $x \in \mathbb{R}^{n_q}$.
		
		\begin{Lemma}[Stability and Gramians]
		\label{lemma:stab}
		$H$ is quadratically stable iff there exist generalized observability Gramians iff there exist generalized controllability Gramians. 
		\end{Lemma}
		\begin{Lemma}\label{observ:lemma}[Observability Gramian and output energy]
		If $\{\cQ_q\}_{q \in Q}$ are observability Gramians, $h_0=(q_0,x_0)$, $(\bq,\bx)=\xi_{H,h_0}(0,w)$, $(\mathbf{o},\by)=\upsilon_{H,h_0}(0,w)$ (i.e. $\bx,\by$ are the continuous state and output trajectories of $H$ if started from the initial state $h_0$ and fed with the timed sequence $w$ and zero continuous input $u=0$), then 
		\[ \int_0^{\infty} \|\by(s)\|^2_2ds \le x_0^T \cQ_{q_0} x_0. \]
		\end{Lemma}

		\begin{Lemma}\label{control:lemma}[Controllability Gramian and input energy]
		If $\{\cP_q\}_{q \in Q}$ are reachability Gramians, $h_0=(q_0,0)$, $(\bq,\bx)=\xi_{H,h_0}(u,w)$ (i.e. $\bx,\bq$ are the continuous and discrete state trajectories of $H$ if started from the initial state $h_0$ and fed with the timed sequence $w$ and continuous input $u$), then 
		\[ \bx(t) \cP^{-1}_{\bq(t)} \bx(t) \le \int_0^{t} \|u(s)\|^2_2.  \]
		\end{Lemma}

		\noindent
		We can formulate the following balanced model reduction.
		\begin{Procedure}
		\label{balancing:proc}
		\begin{enumerate}
		\item Compute  reachabilility and observability Gramians $\{\cP_q >0\}_{q \in Q}$ and $\{\cQ_q >0\}_{q \in Q}$ which satisfy (\ref{contr:gram:eq}), and, respectively (\ref{obs:gram:eq}). 
		\item Find square factor matrices $\bU_q$ so that $\cP_q = \bU_q \bU_q^T$. Additionally, compute the eigenvalue decomposition of the symmetric matrix $\bU_q^T \cQ_q \bU_q$, as 
		\begin{equation*}
		\bU_q^T \cQ_q \bU_q = \bV_q \bLa_q^2 \bV_q^T,
		\end{equation*}
		where 
		\[ \bLa_q=\mathrm{diag}(\sigma_{q,1},\ldots, \sigma_{q,n_q}), \]
		is a diagonal matrix with the real entries sorted in decreasing order, i.e., $\sigma_{q,1} \ge \sigma_{q,2} \ge \cdots \ge \sigma_{q,n_q}$. 
		\item Construct the transformation matrices $\bS_q \in \mathbb{R}^{n_q \times n_q}$ as follows
		\begin{equation}\label{defSq}
		\bS_q = \bLa_q^{1/2} \bV_q^T \bU_q^{-1}.
		\end{equation}
		Define the matrices (with $ q_1=\delta(q_2,\gamma),~q_2\in Q$)
		\begin{equation}
		\label{eq:bal1}
		\begin{split}
		\bar{A}_q &= \bS_q A_q \bS_q^{-1}, \ \ \bar{B}_q = \bS_q B_q, \ \ \bar{C}_q = C_q \bS_q^{-1}, \\ & \hspace{8mm} \bar{M}_{q_2,\gamma,q_1} = \bS_{q_2} M_{q_2,\gamma, q_1} \bS_{q_1}^{-1}.
		\end{split}
		\end{equation}
		\item
		Choose the truncation orders $0 < r_q \le n_q$ and consider the partitioning 
		\small
		\begin{equation}\label{partition_bal}
		\begin{split}
		& \bar{A}_q = \begin{bmatrix} \bar{A}_q^{11} & \bar{A}_q^{12} \\ \bar{A}_q^{21} & \bar{A}_q^{22} 
		\end{bmatrix},  \bar{B}_q = \begin{bmatrix} \bar{B}^1_q  \\ \bar{B}^2_q \end{bmatrix},  \bar{C}_q = \begin{bmatrix} \bar{C}_q^{1} & \bar{C}_q^{2} \end{bmatrix}, ~ r_q < n_q, \\ 
		& \bar{M}_{q_1,\gamma,q_2} = \begin{bmatrix} \bar{M}_{q_1,\gamma,q_2}^{11} & \bar{M}_{q_1,\gamma,q_2}^{12}, \\ \bar{M}_{q_1,\gamma,q_2}^{21} & \bar{M}_{q_1,\gamma,q_2}^{22}
		\end{bmatrix} ~ \mbox{ if } r_{q_1} < n_{q_1}, r_{q_2} < n_{q_2},  \\
		&  \bar{M}_{q_1,\gamma,q_2} = \begin{bmatrix} \bar{M}_{q_1,\gamma,q_2}^{11} & \bar{M}_{q_1,\gamma,q_2}^{12} 
		\end{bmatrix} ~ \mbox{ if } r_{q_1} = n_{q_1}, r_{q_2} < n_{q_2}, \\
		& \bar{M}_{q_1,\gamma,q_2} = \begin{bmatrix} \bar{M}_{q_1,\gamma,q_2}^{11} \\ \bar{M}_{q_1,\gamma,q_2}^{21} 
		\end{bmatrix} ~ \mbox{ if } r_{q_1} < n_{q_1}, r_{q_2} = n_{q_2},
		\end{split}
		\end{equation}
		\normalsize
		where $\bar{A}_q^{11} \in \mathbb{R}^{r_q \times r_q}, \ \bar{M}_{q_1,\gamma,q_2}^{11} \in \mathbb{R}^{r_{q_1} \times r_{q_2}},  \ \bar{B}_q^{1} \in \mathbb{R}^{r_q\times m}, \ \text{and}$ $\bar{C}_q^{1} \in \mathbb{R}^{p \times r_q}$.
		\item
		Define the reduced model
		\begin{align*}
		 \hat{H}=(Q,\Gamma,O,\delta,\lambda, \{r_q,\hat{A}_q,\hat{B}_q,\hat{C}_q \}_{q \in Q},\\
		 \{ \hat{M}_{q_{1},\gamma,q_{2}} \}_{q_{2} \in Q, \gamma \in \Gamma, q_1=\delta(q_2,\gamma)}, (q_0,0)),
		\end{align*}
		 where
		\begin{equation}\label{reduced_bal}
		\begin{split}
		& \hat{A}_q = \bar{A}_q^{11}, \ \ \hat{B}_q = \bar{B}_q^{1}, \ \ \hat{C}_q = \bar{C}_q^{1}, ~ \mbox{ if } r_q \leq n_q, \\ 
		& \hat{M}_{q_1,\gamma,q_2} = \bar{M}_{q_1,\gamma,q_2}^{11}, \mbox{ if }  r_{q_1} < n_{q_1} \mbox{ or }  ~ r_{q_2} < n_{q_2}, \\
		& \hat{A}_q = \bar{A}_q, \ \ \hat{B}_q = \bar{B}_q, \ \ \hat{C}_q = \bar{C}_q,  ~ \mbox{ if } ~ r_q=n_q,\\
		& \hat{M}_{q_1,\gamma,q_2} = \bar{M}_{q_1,\gamma,q_2}, \mbox{ if } r_{q_1} = n_{q_1} \mbox{ and } r_{q_2} = n_{q_2}.  \\
		\end{split}
		\end{equation}  
		\end{enumerate}
		\end{Procedure}
		\begin{Lemma}[Balanced realization]
		\label{balanced:lemma1}
		Consider the \AHLS\ $\bar{H}=(Q,\Gamma,O,\delta,\lambda,$ $\{r_q,\bar{A}_q,\bar{B}_q,\bar{C}_q \}_{q \in Q}, \{ \bar{M}_{q_{1},\gamma,q_{2}} \}_{q_{2} \in Q,\gamma \in \Gamma, q_1=\delta(q_2,\gamma)},$ $(q_0,0))$. Then $\{\bLa_q\}_{q \in Q}$ are both generalized reachability and observability Gramians of $\bar{H}$.
		\end{Lemma}
		
		In the sequel, we will say that an \AHLS\ is \emph{balanced}, if it has generalized reachability Gramians $\{\cP_q\}_{q \in Q}$, generalized observability Gramians 
		$\{\cQ_q\}_{q \in Q}$, and for all $q \in Q$, the matrices $\cQ_q$ and $\cP_q$ are equal and are diagonal. Lemma \ref{balanced:lemma1} says that $\bar{H}$ is balanced.
		In fact, more is true.
		\begin{Lemma}[Preservation of balancing and stability]
		\label{balanced:lemma}
		The reduced order model $\hat{H}$ is balanced, its generalized observability and reachability Gramians are $\{\hat{\bLa}_q\}_{q \in Q}$, $\hat{\bLa}_q=\mathrm{diag}(\sigma_{q,1},\ldots,\sigma_{q,r_q})$. In particular, $\hat{H}$ is quadratically stable. 
		\end{Lemma}

		\begin{Theorem}[Error bound]
		\label{theo:bound}
		For any $(\bu,w) \in \HYBINP$, consider the outputs $(\mathbf{o},\by)=\upsilon_{H}(\bu,w)$ and $(\hat{\mathbf{o}},\hat{\by})=\upsilon_{\hat{H}}(\bu,w)$ generated
		by $H$ and $\hat{H}$ respectively under the input $u$ and timed event sequence $w$ from the corresponding initial state. Then $\hat{\mathbf{o}}=\mathbf{o}$, and 
		\[ \|\by-\hat{\by}\|_{L_2} \le 2(\sum_{q \in Q}\sum_{i=1}^{n_q-r_q} \sigma_{q,r_q+i})\|\bu\|_{L_2}. \]
		\end{Theorem}
		First we prove Theorem \ref{theo:bound} for the case when $n_q-r_q \le 1$ for all $q \in Q$. 
		More precisely, for each $q \in Q$, consider the decomposition 
		\begin{equation}
		\bLa_q = \left[ \begin{array}{cc} \hat{\bLa}_q & 0 \\ 0 & \beta_q
		\end{array} \right],\ \  \beta_q \in \mathbb{R}.
		\end{equation}
		Define $\beta=\min_{q \in Q} \beta_q$ and for each $q \in Q$, define
		\[ r_q = \left\{\begin{array}{rl} n_q-1 & \mbox{ if } \beta_q = \beta, \\
		n_q & \mbox{ otherwise }
		\end{array}\right. .
		\]
		Consider the reduced order model $\hat{H}$ from Procedure \ref{balancing:proc} for this choice of $r_q$. 
		\begin{Theorem}[One step error bound]
		\label{theo:bound:step1}
		For any $(\bu,w) \in \HYBINP$, consider the outputs $(\mathbf{o},\by)=\upsilon_{H}(\bu,w)$ and $(\hat{\mathbf{o}},\hat{\by})=\upsilon_{\hat{H}}(\bu,w)$ generated
		by $H$ and $\hat{H}$ respectively under the input $u$ and timed event sequence $w$ from the corresponding initial state. Then $\hat{\mathbf{o}}=\mathbf{o}$, and 
		\[ \|\by-\hat{\by}\|_{L_2} \le 2\beta\|\bu\|_{L_2}. \]
		\end{Theorem}
		Theorem \ref{theo:bound} follows by repeated application of Theorem \ref{theo:bound:step1}. 
		The proof of Theorem \ref{theo:bound:step1} is done via a sequence of lemmas.  In order to state these lemmas, we introduce the following notation.
		Consider the balanced \AHLS $\bar{H}$ from Lemma \ref{balanced:lemma1}. Note that the \SAHLS\ $\bar{H}$ and $H$ are isomorphic, and hence they have the same input-output map.
		Consider now the state trajectory $(\bq,\bar{\bx})=\xi_{\bar{H},h_0}(\bu,w)$  of $\bar{H}$ and the state trajectory $(\hat{\bq},\hat{\bx})=\xi_{\hat{H},\hat{h}_0}(\bu,w)$, 
		$\hat{h}_0=(q_0,0)$ is the initial state of $\hat{H}$. It is easy to see that $\bq=\hat{\bq}$. 
		
		For any $t \in \mathbb{R}_{+}$ such that $r_{\bq(t)}=n_{\bq(t)}-1$, consider the partitioning
		\[ \bar{\bx}(t) = \left[ \begin{array}{c} \bar{\bx}_1(t) \\ \bar{\bx}_2(t) \end{array} \right], \]
		with $\bar{\bx}_1(t) \in \mathbb{R}^{r_{q_i}}, \ \bar{\bx}_2(t) \in \mathbb{R}$.
		Define the functions
		\begin{equation}\label{xoxc}
		\begin{split}
		& \bx_o(t) = \left\{\begin{array}{rl} \left[ \begin{array}{c} \bar{\bx}_1(t) - \hat{\bx}(t) \\ \bar{\bx}_2(t) \end{array} \right],  & r_{\bq(t)}=n_{\bq(t)}-1 \\ 
		\bar{\bx}(t)-\hat{\bx}(t) & \mbox{ otherwise }
		\end{array}\right. ,\\
		&  \bx_c(t) = \left\{\begin{array}{rl} \left[ \begin{array}{c} \bar{\bx}_1(t) + \hat{\bx}(t) \\ \bar{\bx}_2(t) \end{array} \right], & r_{\bq(t)}=n_{\bq(t)}-1 \\
		\bar{\bx}(t)+\hat{\bx}(t) & \mbox{ otherwise }
		\end{array}\right. .\\
		\end{split}
		\end{equation}
		Note that the following holds:
		\begin{equation*}
		\by(t) - \hat{\by}(t) = C_{\bq(t)} \bx_o(t).
		\end{equation*}
		Define the function 
		\small
		\begin{equation}\label{defVxoxc}
		V(\bx_o(t), \bx_c(t)) = \bx_o(t)^T \bLa_{\bq(t)} \bx_o(t)+ \beta^2  \bx_c(t)^T \bLa_{\bq(t)}^{-1} \bx_c(t).
		\end{equation}
		\normalsize
		\begin{Lemma} 
		\label{error:proof:lemma1}
		The temporal derivative of the function V, as defined in (\ref{defVxoxc}), satisfies
		\begin{equation}\label{derVineq}
		\frac{\partial V(\bx_o(t), \bx_c(t))}{\partial t} \leqslant  4 \beta^2 \Vert \bu(t) \Vert_2^2- \Vert \by(t)-\hat{\by}(t) \Vert_2^2,
		\end{equation}
		for all $t \in [T_{i-1}, T_i)$.
		\end{Lemma}
		\begin{IEEEproof}[Proof of Lemma \ref{error:proof:lemma1}]
		Note that 
		\begin{align}
		\bar{A}_q \bLa_q +\bLa_q \bar{A}^T_q+\bar{B}_q \bar{B}_q^T < 0, 
		\label{Lyap_reach_bal} \\
		\bar{A}_q^T \bLa_q +\bLa_q \bar{A}_q+ \bar{C}_q^T \bar{C}_q < 0. ~ 
		\label{Lyap_obser_bal}
		\end{align}
		Two cases have to be distinguished.

		The first one is when $r_{q_i}=n_{q_i}$, i.e., in the discrete mode $q_i$ no truncation takes place. In that case, notice that
		\begin{align}
		& \dot{\bx}_o(t) = \bar{A}_{q_i} \bx_o (t), 
		& \dot{\bx}_c(t) = \bar{A}_{q_i} \bx_c (t) +2 \bar{B}_{q_i}^2 \bu(t). 
		\end{align}
		We observe that  $\frac{d }{d t} \bx_o(t)^T \bLa_{q_i} \bx_o(t) = 2 (\bar{A}_{q_i} \bx_o (t))^T \bLa_{q_i} \bx_o(t) \le -\bx^T_o(t) \bar{C}^T_{q_i}\bar{C}_{q_i}\bx_o(t)=-\|\by(t)-\hat{\by}(t)\|^2_2$
		due to \eqref{Lyap_obser_bal} and Remark \ref{LMI:rem1}. By Remark \ref{LMI:rem1} and \eqref{Lyap_reach_bal}, 
		$\frac{d }{d t} \bx_c(t)^T \bLa_{q_i}^{-1} \bx_c(t) = 2 (\bar{A}_{q_i} \bx_c (t) + 2\bar{B}_{q_i}u(t))^T \bLa_{q_i}^{-1} \bx_c(t) \le -4\|\bu(t)\|^2_2$. Hence,
		the claim of the lemma is satisfied. 
		
		Assume now that $r_{q_i}=n_{q_i}-1$. Then $\beta_{q_i}=\beta$ and the following holds:
		\small
		\begin{align}
		\dot{\bx}_o(t) &= \bar{A}_{q_i} \bx_o (t) + \left[ \begin{array}{c} \bfz \\ \bar{B}_{q_i}^2(t) \end{array} \right] \bu(t) + \left[ \begin{array}{c} \bfz \\ \bar{A}_{q_i}^{21}(t) \end{array} \right] \hat{\bx}(t), \label{xo_der} \\
		\dot{\bx}_c(t) &= \bar{A}_{q_i} \bx_c (t) +2 \bar{B}_{q_i} \bu(t) - \left[ \begin{array}{c} \bfz \\ \bar{B}_{q_i}^2(t) \end{array} \right] \bu(t) -  \left[ \begin{array}{c} \bfz \\ \bar{A}_{q_i}^{21}(t) \end{array} \right] \hat{\bx}(t). \label{xc_der}
		\end{align}
		\normalsize
		By using (\ref{xo_der}), \eqref{Lyap_obser_bal}, \eqref{obs:gram:eq2} and Remark \ref{LMI:rem1}, it follows that
		\begin{align}\label{derVxo}
	    \begin{split}
		& \frac{d }{d t} \bx_o(t)^T \bLa_{q_i} \bx_o(t) = 2 \bx_o^T (t)\bar{A}_{q_i}^T \bLa_{q_i} \bx_o(t)\\ &\qquad+ 2 \left( \left[ \begin{array}{c} \bfz \\ \bar{B}_{q_i}^2 \bu(t)+\bar{A}_{q_i}^{21} \hat{\bx}(t) \end{array} \right]^T \bLa_{q_i} \bx_o(t) \right)  \\
		& \leqslant \Vert \bar{C}_{q_i} \bx_o(t)  \Vert_2^2 + 2 \alpha_o  = - \Vert \by(t) - \hat{\by}(t) \Vert_2^2  +  2 \alpha_o,
	    \end{split}
		\end{align}
		where  
		\begin{align} \label{alphao}
	    \begin{split}
		\alpha_o  &= \left[ \begin{array}{c} \bfz \\ \bar{B}_{q_i}^2 \bu(t)+\bar{A}_{q_i}^{21} \hat{\bx}(t) \end{array} \right]^T  \left[ \begin{array}{cc} \hat{\bLa}_{q_i} & \bfz \\ \bfz & \beta_{q_i} \end{array} \right]  \left[ \begin{array}{c} \bar{\bx}_1(t)- \hat{\bx}(t) \\ \bar{\bx}_2(t) \end{array} \right]  \\
		&=   \beta_{q_i} \big{(} \bar{B}_{q_i}^2 \bu(t)+\bar{A}_{q_i}^{21} \hat{\bx}(t) \big{)}^T \bar{\bx}_2(t).
	    \end{split}
		\end{align}
		Similarly, by using (\ref{xc_der}), \eqref{contr:gram:eq2} from Remark \ref{LMI:rem1} and \eqref{Lyap_reach_bal}, we show that
		\begin{align}\label{derVxc}
	    \begin{split}
		& \frac{d }{d t} \bx_c(t)^T \bLa_{q_i}^{-1} \bx_c(t) = 2 \big{(} \bar{A}_{q_i}  \bx_c(t) + \bar{B}_{q_i} 2 \bu(t)  \big{)}^T \bLa_{q_i}^{-1} \bx_c(t) \\ &-  2 \left( \left[ \begin{array}{c} \bfz \\ \bar{B}_{q_i}^2 \bu(t)+\bar{A}_{q_i}^{21} \hat{\bx}(t) \end{array} \right]^T \bLa_{q_i}^{-1} \bx_c(t) \right) 
		\leqslant 4 \Vert \bu(t) \Vert_2^2- 2 \alpha_c ,
		\end{split}
		\end{align}
		where 
		\begin{align} \label{alphac}
		\begin{split}
		\alpha_c  &= \left[ \begin{array}{c} \bfz \\ \bar{B}_{q_i}^2 \bu(t)+\bar{A}_{q_i}^{21} \hat{\bx}(t) \end{array} \right]^T \left[ \begin{array}{cc} \hat{\bLa}_{q_i}^{-1} & \bfz \\ \bfz & \beta_{q_i}^{-1}  \end{array} \right]  \left[ \begin{array}{c} \bar{\bx}_1(t)+ \hat{\bx}(t) \\ \bar{\bx}_2(t) \end{array} \right]   \\ &= 
		\beta_{q_i}^{-1} \big{(} \bar{B}_{q_i}^2 \bu(t)+\bar{A}_{q_i}^{21} \bar{\bx}(t) \big{)}^T \bar{\bx}_2(t).
		\end{split}
		\end{align}
		From (\ref{alphao}) and (\ref{alphac}) and $\beta=\beta_{q_i}$, observe that $\alpha_o = \beta^2 \alpha_c$. Hence, by adding the inequality in (\ref{derVxo}) with the one in (\ref{derVxc}) multiplied by $\beta^2=\beta_{q_i}^2$, it follows that
		\begin{multline*}
		\frac{d }{d t} \bx_o(t)^T \bLa_{q_i} \bx_o(t) + \beta^2  \frac{d }{d t} \bx_c(t)^T \bLa_{q_i}^{-1} \bx_c(t) \\  \leqslant  - \Vert \by(t) -  \hat{\by}(t) \Vert_2^2 + 4 \beta_{q_i}^2 \Vert \bu(t) \Vert_2^2,
		\end{multline*}
		and by using the definition of $V$ in (\ref{defVxoxc}), it automatically proves the result in (\ref{derVineq}).
		\end{IEEEproof}
		\begin{Lemma}
		\label{error:proof:lemma2}
		For all $i \in \mathbb{N}$, 
		\begin{equation}
		\label{jumpVineq}
		V(\bx(T_{i+1}),\hat{\bx}(T_{i+1})) \le V(\bx(T_{i+1}^{-}),\hat{\bx}(T_{i+1}^{-})),
		\end{equation}
		where $\bx(T_{i+1}^{-})=\lim_{t \uparrow T_{i+1}} \bx(t)$, and $\hat{\bx}(T_{i+1}^{-})=\lim_{t \uparrow T_{i+1}} \hat{\bx}(t)$. 
		\end{Lemma}	
		
		\begin{IEEEproof}[Proof of Lemma \ref{error:proof:lemma2}]
		Note that $q_i=\bq(t)$ for all $t \in [T_i,T_{i+1})$ and that $\delta(q_i,\gamma_{i+1})=q_{i+1}$. Moreover, by virtue of $\{\bLa_q\}_{q \in Q}$ being 
		generalized observability and reachability Gramians for $\bar{H}$, and Remark \ref{LMI:rem1}, the following holds
		\begin{align}
		\bar{M}_{q_{i+1},\gamma_{i+1},q_i}^T \bLa_{q_{i+1}}^{-1}\bar{M}_{q_{i+1},\gamma_{i+1},q_i}  < \bLa_{q_i}^{-1}, 
		\label{Lyap_reach_bal:reset} \\
		\bar{M}_{q_{i+1},\gamma_{i+1},q_i}^T \bLa_{q_{i+1}}\bar{M}_{q_{i+1},\gamma_{i+1},q_i}  < \bLa_{q_i}. 
		\label{Lyap_obser_bal:reset}
		\end{align}
		
		In order to prove \eqref{jumpVineq}, the following cases have to be distinguished. 
		
		Assume that $r_{q_{i+1}}=n_{q_i+1}$, i.e., no truncation takes place in mode $q_{i+1}$. In this case,
		$\bx(T_{i+1}) =\bar{M}_{q_{i+1},\gamma_{i+1},q_i}\bx(T_{i+1}^{-})$, and 
		\begin{equation}
		\label{case1:jump:eq-3}
		 \begin{split}
		  \hat{\bx}(T_{i+1})= \bar{M}_{q_{i+1},\gamma_{i+1},q_i}^{11}  \hat{\bx}(T_{i+1}^{-}) = 
		\bar{M}_{q_{i+1},\gamma_{i+1},q_i} \begin{bmatrix} \hat{\bx}(T_{i+1}^{-}) \\ 0 \end{bmatrix} ,
		\end{split}
		\end{equation}
		if $r_{q_i}=n_{q_i}-1$, and 
		\begin{equation}
		\label{case1:jump:eq-3.1}
		\hat{\bx}(T_{i+1})= \bar{M}_{q_{i+1},\gamma_{i+1},q_i} \hat{\bx}(T_{i+1}^{-}),
		\end{equation}
		if $r_{q_i}=n_{q_i}$. 
		Notice that if $r_{q_i}=n_{q_i}$, then
		\begin{equation}
		\label{case1:jump:eq-2}
		\begin{split}
		& \bx_{c}(T_{i+1})=\bx(T_{i+1})+\hat{\bx}(T_{i+1}), \\ & \bx_{o}(T_{i+1})=\bx(T_{i+1})-\hat{\bx}(T_{i+1}),  \\
		& \bx_{c}(T_{i+1}^{-})=\bx(T_{i+1}^{-})+\hat{\bx}(T_{i+1}^{-}),\\ 
		& \bx_{o}(T_{i+1}^{-})=\bx(T_{i+1}^{-})-\hat{\bx}(T_{i+1}^{-}).
		\end{split}
		\end{equation}
		Similarly, if $r_{q_i} = n_{q_i}-1$, then
		\begin{equation}
		\label{case1:jump:eq-1}    
		\begin{split}
		& \bx_{c}(T_{i+1}^{-})=\bx(T_{i+1}^{-})+\begin{bmatrix} \hat{\bx}(T_{i+1}^{-}) \\ 0 \end{bmatrix},  \\
		&  \bx_{o}(T_{i+1}^{-})=\bx(T_{i+1}^{-})-\begin{bmatrix} \hat{\bx}(T_{i+1}^{-}) \\ 0 \end{bmatrix}.
		\end{split}
		\end{equation}
		From \eqref{case1:jump:eq-3}-\eqref{case1:jump:eq-1}, 
		it follows that
		\begin{equation}
		\label{case1:jump:eq}
		 \begin{split}
		& \bx_c(T_{i+1})=\bar{M}_{q_{i+1},\gamma_{i+1},q_i}\bx_c(T_{i+1}^{-}),  \\
		& \bx_o(T_{i+1})=\bar{M}_{q_{i+1},\gamma_{i+1},q_i}\bx_o(T_{i+1}^{-}).
		\end{split}
		\end{equation}     
		From \eqref{case1:jump:eq} it then follows that
		\small
		\begin{multline}
		\label{case1:jump:eq2}
		        V(\bx(T_{i+1}),\hat{\bx}(T_{i+1})) =  
		\bx_{o}^T(T_{i+1}^{-}) \bar{M}_{q_{i+1},\gamma_{i+1},q_i}^T \bLa_{q_{i+1}}\bar{M}_{q_{i+1},\gamma_{i+1},q_i}\bx_{o}(T_{i+1}^{-}) \\ + 
		        \beta^{2} \bx_c^T(T_{i+1}^{-})\bar{M}_{q_{i+1},\gamma_{i+1},q_i}^T \bLa_{q_{i+1}}^{-1}\bar{M}_{q_{i+1},\gamma_{i+1},q_i}\bx_c(T_{i+1}^{-}). 
		\end{multline}
		\normalsize 
		From \eqref{Lyap_obser_bal:reset}-\eqref{Lyap_reach_bal:reset} it follows that
		\small
		\[
		\begin{split}
		& \bx_{o}^T(T_{i+1}^{-}) \bar{M}_{q_{i+1},\gamma_{i+1},q_i}^T \bLa_{q_{i+1}}\bar{M}_{q_{i+1},\gamma_{i+1},q_i}\bx_{o}(T_{i+1}^{-}) \le \bx^T_{o}(T_{i+1}^{-})\bLa_{q_i}\bx_{o}(T_{i+1}^{-}), \\
		& \bx_{c}^T(T_{i+1}^{-}) \bar{M}_{q_{i+1},\gamma_{i+1},q_i}^T \bLa_{q_{i+1}}^{-1}\bar{M}_{q_{i+1},\gamma_{i+1},q_i}\bx_{c}(T_{i+1}^{-}) \le \bx^T_{c}(T_{i+1}^{-})\bLa_{q_i}^{-1}\bx_{c}(T_{i+1}^{-}).
		\end{split}
		\]
		\normalsize
		Hence, from \eqref{case1:jump:eq2}, it follows that 
		\begin{multline*}
		V(\bx(T_{i+1}),\hat{\bx}(T_{i+1})) \le  \bx^T_{o}(T_{i+1}^{-})\bLa_{q_i}\bx_{o}(T_{i+1}^{-}) \\ + \beta^2 \bx^T_{c}(T_{i+1}^{-})\bLa_{q_i}^{-1}\bx_{c}(T_{i+1}^{-}) 
		= V(\bx(T_{i+1}^{-}),\hat{\bx}(T_{i+1}^{-})), 
		\end{multline*}
		i.e.,  \eqref{jumpVineq} holds. 
		
		Consider now the case when $r_{q_{i+1}}=n_{q_{i+1}}-1$, i.e., in mode $q_{i+1}$ truncation takes place. In this case, 
		$\bx(T_{i+1})=\bar{M}_{q_{i+1},\gamma_{i+1},q_i}\bx(T_{i+1}^{-})$, and 
		\begin{equation}
		\label{case2:jump:eq-5}
		\begin{split}
		\hat{\bx}(T_{i+1})&=
		\bar{M}_{q_{i+1},\gamma_{i+1},q_i}^{11}  \hat{\bx}(T_{i+1}^{-})\\
		&= \bar{M}_{q_{i+1},\gamma_{i+1},q_i} \begin{bmatrix} \hat{\bx}(T_{i+1}^{-})   \\ 0 \end{bmatrix}  - \begin{bmatrix} 0 \\ \bar{M}_{q_{i+1},\gamma_{i+1},q_i}^{21}	\end{bmatrix} \hat{\bx}(T_{i+1}^{-}),
		\end{split}
		\end{equation}
		if $r_{q_i}=n_{q_i}-1$, and 
     \begin{equation}
       \label{case2:jump:eq-4}
     \begin{split}
     \hat{\bx}(T_{i+1})&=
          \bar{M}_{q_{i+1},\gamma_{i+1},q_i}^{11}  \hat{\bx}(T_{i+1}^{-})\\
          &= \bar{M}_{q_{i+1},\gamma_{i+1},q_i} \hat{\bx}(T_{i+1}^{-})  - \begin{bmatrix} 0 \\ \bar{M}_{q_{i+1},\gamma_{i+1},q_i}^{21}	\end{bmatrix} \hat{\bx}(T_{i+1}^{-}), 
         \end{split} 
      \end{equation}
      if $r_q=n_q$. 
		Notice that
		\begin{equation}
		\label{case2:jump:eq-3}
		\begin{split}
		& \bx_{c}(T_{i+1})=\bx(T_{i+1})+\begin{bmatrix} \hat{\bx}(T_{i+1}) \\ 0 \end{bmatrix}, \\ 
		& \bx_{o}(T_{i+1})=\bx(T_{i+1})-\begin{bmatrix} \hat{\bx}(T_{i+1}) \\ 0 \end{bmatrix},  
	   \end{split}
	   \end{equation}
	   and if $r_q=n_q$, then
	   \begin{equation}
	   \label{case2:jump:eq-2}
	     \begin{split}
		& \bx_{c}(T_{i+1}^{-})=\bx(T_{i+1}^{-})+\hat{\bx}(T_{i+1}^{-}),  \\ & \bx_{o}(T_{i+1}^{-})=\bx(T_{i+1}^{-})-\hat{\bx}(T_{i+1}^{-}), \\
	    \end{split}
	    \end{equation}
	    and for $r_q=n_q-1$,
	    \begin{equation}
	   \label{case2:jump:eq-1}
	    \begin{split}
		&\bx_{c}(T_{i+1}^{-})=\bx(T_{i+1}^{-})+\begin{bmatrix}  \hat{\bx}(T_{i+1}^{-}) \\ 0 \end{bmatrix}, \\ & \bx_{o}(T_{i+1}^{-}) 
		 =\bx(T_{i+1}^{-})-\begin{bmatrix} \hat{\bx}(T_{i+1}^{-}) \\ 0 \end{bmatrix}. 
		\end{split}
		\end{equation}
		From \eqref{case2:jump:eq-5}-\eqref{case2:jump:eq-1} it then follows that
		\begin{equation}
		\label{case2:jump:eq2}
		\begin{split}
		\bx_c(T_{i+1})=\bar{M}_{q_{i+1},\gamma_{i+1},q_i}\bx_c(T_{i+1}^{-}) - \begin{bmatrix} 0 \\ \bar{M}^{21}_{q_{i+1},\gamma_{i+1},q_i} \end{bmatrix}\hat{\bx}(T_{i+1}^{-}), \\
		\bx_o(T_{i+1})=\bar{M}_{q_{i+1},\gamma_{i+1},q_i}\bx_o(T_{i+1}^{-}) + \begin{bmatrix} 0 \\ \bar{M}_{q_{i+1},\gamma_{i+1},q_i}^{21} \end{bmatrix}\hat{\bx}(T_{i+1}^{-}).  
		\end{split}
		\end{equation}     
		From \eqref{case2:jump:eq2} it then follows that
		\begin{equation}
		\label{case2:obs:eq1}
		\begin{split}
		& \bx_{o}^T(T_{i+1})\bLa_{q_{i+1}}\bx_{o}(T_{i+1}) =  \\
		& \bx_{o}^T(T_{i+1}^{-}) \bar{M}_{q_{i+1},\gamma_{i+1},q_i}^T \bLa_{q_{i+1}}\bar{M}_{q_{i+1},\gamma_{i+1},q_i}
		\bx_{o}(T_{i+1}^{-})+          \\ 
		& 2 \bx_{o}^T(T_{i+1}^{-1})\bar{M}_{q_{i+1},\gamma_{i+1},q_i}^T\bLa_{q_{i+1}}\begin{bmatrix} 0 \\ \bar{M}^{21}_{q_{i+1},\gamma_{i+1},q_i} \end{bmatrix}\hat{\bx}(T_{i+1}^{-}) \\
		&+ \left(\begin{bmatrix} 0 \\ \bar{M}^{21}_{q_{i+1},\gamma_{i+1},q_i} \end{bmatrix}\hat{\bx}(T_{i+1}^{-})\right)^T\bLa_{q_{i+1}}\begin{bmatrix} 0 \\ \bar{M}^{21}_{q_{i+1},\gamma_{i+1},q_i} \end{bmatrix}\hat{\bx}(T_{i+1}^{-}).
		\end{split}
		\end{equation}
		Since $\bLa_{q_{i+1}}=\begin{bmatrix} \hat{\bLa}_{q_{i+1}} & 0 \\ 0 & \beta_{q_{i+1}} \end{bmatrix}$, 
		it follows that
		\begin{multline*}
		 \left(\begin{bmatrix} 0 \\ \bar{M}^{21}_{q_{i+1},\gamma_{i+1},q_i} \end{bmatrix}\hat{\bx}(T_{i+1}^{-})\right)^T\bLa_{q_{i+1}}\begin{bmatrix} 0 \\  \bar{M}_{q_{i+1},\gamma_{i+1},q_i}^{21} \end{bmatrix}\hat{\bx}(T_{i+1}^{-}) \\ = \beta_{q_{i+1}} \| \bar{M}_{q_{i+1},\gamma_{i+1},q_i}^{21} \hat{\bx}(T_{i+1}^{-}) \|^2_2.
		\end{multline*}
		Moreover,
		 \begin{multline*}
	    2\bx_{o}^T(T_{i+1}^{-})\bar{M}_{q_{i+1},\gamma_{i+1},q_i}^T\bLa_{q_{i+1}}\begin{bmatrix} 0 \\ \bar{M}^{21}_{q_{i+1},\gamma_{i+1},q_i} \end{bmatrix}\hat{\bx}(T_{i+1}^{-})= \\
		  \gamma_o-2\beta_{q_{i+1}} \| \bar{M}^{21}_{q_{i+1},\gamma_{i+1},q_i} \hat{\bx}(T_{i+1}^{-}) \|^2_2,
		\end{multline*}
		where 
		\[
		\gamma_o= \left\{\begin{array}{rl} 
		& \hspace{-4mm} 2\beta_{q_{i+1}} \left(\bar{M}_{q_{i+1},\gamma_{i+1},q_i}^{21}\bx_1(T_{i+1}^-)+\bar{M}_{q_{i+1},\gamma_{i+1},q_i}^{22}\bx_2(T_{i+1}^{-})\right)^T  \\[1mm] &   \hspace{-4mm}  \times \bar{M}^{21}_{q_{i+1},\gamma_{i+1},q_i}\hat{\bx}(T_{i+1}^{-})   \mbox{ if } r_{q_{i}}=n_{q_i}-1 \\[1mm]
		& \hspace{-4mm} 2\beta_{q_{i+1}}  \left(\bar{M}_{q_{i+1},\gamma_{i+1},q_i}^{21} \bx(T_{i+1}^-)\right)^T  \\[1mm] & \hspace{-4mm} \times \bar{M}^{21}_{q_{i+1},\gamma_{i+1},q_i}\hat{\bx}(T_{i+1}^{-})   \mbox{ if } r_{q_{i}}=n_{q_i} \end{array}\right. . 
		\]
		Hence, it follows that
		\begin{equation}
		\label{case2:obs:eq2}
		\begin{split}
		&\bx_{o}^T(T_{i+1})\bLa_{q_{i+1}}\bx_{o}(T_{i+1})\\
		&= \bx_{o}^T(T_{i+1}^{-}) \bar{M}_{q_{i+1},\gamma_{i+1},q_i}^T \bLa_{q_{i+1}}\bar{M}_{q_{i+1},\gamma_{i+1},q_i} \bx_{o}(T_{i+1}^{-}) \\
		 &\quad + \gamma_o-\beta_{q_{i+1}}\| \bar{M}^{21}_{q_{i+1},\gamma_{i+1},q_i} \hat{\bx}(T_{i+1}^{-}) \|^2_2.
		\end{split}
		\end{equation}
		With a similar reasoning, 
		\begin{equation}
		\label{case2:reach:eq1}
		\begin{split}
		& \bx_{c}^T(T_{i+1})\bLa^{-1}_{q_{i+1}}\bx_{c}(T_{i+1})\\
		&= \bx_{c}^T(T_{i+1}^{-}) \bar{M}_{q_{i+1},\gamma_{i+1},q_i}^T \bLa_{q_{i+1}}^{-1}   \bar{M}_{q_{i+1},\gamma_{i+1},q_i}
		 \bx_{c}(T_{i+1}^{-})     \\ &~~ -2\bx_{c}^T(T_{i+1}^{-1})\bar{M}_{q_{i+1},\gamma_{i+1},q_i}^T\bLa^{-1}_{q_{i+1}}\begin{bmatrix} 0 \\ \bar{M}^{21}_{q_{i+1},\gamma_{i+1},q_i} \end{bmatrix}\hat{\bx}(T_{i+1}^{-}) \\
		&~~+ (\begin{bmatrix} 0 \\ \bar{M}^{21}_{q_{i+1},\gamma_{i+1},q_i} \end{bmatrix}\hat{\bx}(T_{i+1}^{-}))^T\bLa_{q_{i+1}}^{-1}\begin{bmatrix} 0 \\ \bar{M}^{21}_{q_{i+1},\gamma_{i+1},q_i} \end{bmatrix}\hat{\bx}(T_{i+1}^{-}).
		\end{split}
		\end{equation}
		Since $\bLa_{q_{i+1}}^{-1}=\begin{bmatrix} \hat{\bLa}_{q_{i+1}}^{-1} & 0 \\ 0 & \beta_{q_{i+1}}^{-1} \end{bmatrix}$,  we can again write that 
		\begin{multline*}
		\left(\begin{bmatrix} 0 \\ \bar{M}^{21}_{q_{i+1},\gamma_{i+1},q_i} \end{bmatrix}\hat{\bx}(T_{i+1}^{-})\right)^T\bLa_{q_{i+1}}^{-1}\begin{bmatrix} 0 \\ \bar{M}^{21}_{q_{i+1},\gamma_{i+1},q_i} \end{bmatrix}\hat{\bx}(T_{i+1}^{-}) \\ = \beta^{-1}_{q_{i+1}} \| \bar{M}^{21}_{q_{i+1},\gamma_{i+1},q_i} \hat{\bx}(T_{i+1}^{-}) \|^2_2, 
		   \end{multline*}
		 and 
		\begin{multline*}
		2\bx_{c}^T(T_{i+1}^{-})\bar{M}_{q_{i+1},\gamma_{i+1},q_i}^T\bLa_{q_{i+1}}^{-1}\begin{bmatrix} 0 \\ \bar{M}^{21}_{q_{i+1},\gamma_{i+1},q_i} \end{bmatrix}\hat{\bx}(T_{i+1}^{-})\\
		= \gamma_c+2\beta_{q_{i+1}}^{-1} \| \bar{M}^{21}_{q_{i+1},\gamma_{i+1},q_i} \hat{\bx}(T_{i+1}^{-}) \|^2_2,
		\end{multline*}
		where 
		\[
		\gamma_c= \left\{\begin{array}{rl} 
		& \hspace{-4mm} 2\beta_{q_{i+1}}^{-1} \left(\bar{M}_{q_{i+1},\gamma_{i+1},q_i}^{21}\bx_1(T_{i+1}^-)+\bar{M}_{q_{i+1},\gamma_{i+1},q_i}^{22}\bx_2(T_{i+1}^{-})\right)^T  \\[1mm] & \hspace{-4mm} \times \bar{M}^{21}_{q_{i+1},\gamma_{i+1},q_i}\hat{\bx}(T_{i+1}^{-})   \mbox{ if } r_{q_{i}}=n_{q_i}-1 \\[1mm]
		& \hspace{-4mm} 2\beta_{q_{i+1}}^{-1}  \left(\bar{M}_{q_{i+1},\gamma_{i+1},q_i}^{21} \bx(T_{i+1}^-)\right)^T  \\[1mm]
		& \hspace{-4mm} \times \bar{M}^{21}_{q_{i+1},\gamma_{i+1},q_i}\hat{\bx}(T_{i+1}^{-})   \mbox{ if } r_{q_{i}}=n_{q_i} \end{array}\right. ,
		\]
		and hence
		\begin{equation}
		\label{case2:reach:eq2}
		\begin{split}
		& \bx_{c}^T(T_{i+1})\bLa_{q_{i+1}}^{-1}\bx_{c}(T_{i+1})\\ 
		&= \bx_{c}^T(T_{i+1}^{-}) \bar{M}_{q_{i+1},\gamma_{i+1},q_i}^T \bLa_{q_{i+1}}^{-1}\bar{M}_{q_{i+1},\gamma_{i+1},q_i} 
		 \bx_{c}(T_{i+1}^{-})\\
		 &\quad -\gamma_c-\beta_{q_{i+1}}^{-1}\| \bar{M}^{21}_{q_{i+1},\gamma_{i+1},q_i} \hat{\bx}(T_{i+1}^{-}) \|^2_2.
		\end{split}
		\end{equation}
		Note that $\beta=\beta_{q_{i+1}}$ since it was assumed that $r_{q_{i+1}}=n_{q_{i+1}}-1$. Moreover, notice that $\beta^2_{q_{i+1}}\gamma_c=\gamma_o$, hence
		by using \eqref{case2:obs:eq2} and \eqref{case2:reach:eq2}
		\begin{equation*}
		\begin{split}
		& V(\bx(T_{i+1}),\hat{\bx}(T_{i+1}))\\ 
		&=  \bx_{o}^T(T_{i+1})\bLa_{q_{i+1}}\bx_{o}(T_{i+1})
		+ \beta^2 \bx_{c}^T(T_{i+1})  \bLa^{-1}_{q_{i+1}}\bx_{c}(T_{i+1}) \\
		&=  \bx_{o}^T(T_{i+1}^{-}) \bar{M}_{q_{i+1},\gamma_{i+1},q_i}^T \bLa_{q_{i+1}}\bar{M}_{q_{i+1},\gamma_{i+1},q_i}\bx_{o}(T_{i+1}^{-})\\
		&\quad+\gamma_o-\beta_{q_{i+1}}\| \bar{M}^{21}_{q_{i+1},\gamma_{i+1},q_i} \hat{\bx}(T_{i+1}^{-}) \|^2_2 \\
		&\quad +\beta^2_{q_{i+1}} \bx_{c}^T(T_{i+1}^{-}) \bar{M}_{q_{i+1},\gamma_{i+1},q_i}^T
		\bLa_{q_{i+1}}^{-1}\bar{M}_{q_{i+1},\gamma_{i+1},q_i}\bx_{c}(T_{i+1}^{-}) \\
		&\quad -\beta^2_{q_{i+1}}\gamma_c-\beta^2_{q_{i+1}} \beta_{q_{i+1}}^{-1} \|\bar{M}^{21}_{q_{i+1},\gamma_{i+1},q_i} \hat{\bx}(T_{i+1}^{-}) \|^2_2, 
		\end{split}
		\]
		and therefore
		\[
		\begin{split}
		& V(\bx(T_{i+1}),\hat{\bx}(T_{i+1})) \\
		& =\bx_{o}^T(T_{i+1}^{-}) \bar{M}_{q_{i+1},\gamma_{i+1},q_i}^T \bLa_{q_{i+1}}\bar{M}_{q_{i+1},\gamma_{i+1},q_i}\bx_{o}(T_{i+1}^{-})\\ 
		&\quad +
		\beta^{2} \bx_c^T(T_{i+1}^{-})\bar{M}_{q_{i+1},\gamma_{i+1},q_i}^T \bLa_{q_{i+1}}^{-1}\bar{M}_{q_{i+1},\gamma_{i+1},q_i}\bx_c(T_{i+1}^{-}) \\ 
		&\quad -2\beta\| \bar{M}^{21}_{q_{i+1},\gamma_{i+1},q_i} \hat{\bx}(T_{i+1}^{-}) \|^2_2. \\
		\end{split}
		\end{equation*} 
		Using that $2\beta^2\|\bar{M}^{21}_{q_{i+1},\gamma_{i+1},q_i} \hat{\bx}(T_{i+1}^{-}) \|^2_2 \ge 0$, it then follows that
		\[
		\begin{split}
		& V(\bx(T_{i+1}),\hat{\bx}(T_{i+1})) \\ 
		&\le \bx_{o}^T(T_{i+1}^{-}) \bar{M}_{q_{i+1},\gamma_{i+1},q_i}^T \bLa_{q_{i+1}}\bar{M}_{q_{i+1},\gamma_{i+1},q_i}\bx_{o}(T_{i+1}^{-})\\
		&\quad+
		\beta^{2} \bx_c^T(T_{i+1}^{-})\bar{M}_{q_{i+1},\gamma_{i+1},q_i}^T \bLa_{q_{i+1}}^{-1}\bar{M}_{q_{i+1},\gamma_{i+1},q_i}\bx_c(T_{i+1}^{-}).
		\end{split}
		\]
		From \eqref{Lyap_obser_bal:reset} and \eqref{Lyap_reach_bal:reset}, it then  follows that
		\[
		\begin{split}
		& V(\bx(T_{i+1}),\hat{\bx}(T_{i+1})) \\ &\le \bx_{o}^T(T_{i+1}^{-}) \bar{M}_{q_{i+1},\gamma_{i+1},q_i}^T \bLa_{q_{i+1}}\bar{M}_{q_{i+1},\gamma_{i+1},q_i} \bx_{o}(T_{i+1}^{-}) \\
		&\quad +\beta^{2} \bx_c^T(T_{i+1}^{-})\bar{M}_{q_{i+1},\gamma_{i+1},q_i}^T \bLa_{q_{i+1}}^{-1}\bar{M}_{q_{i+1},\gamma_{i+1},q_i}\bx_c(T_{i+1}^{-}) \\
		& \le  \bx_{o}^T(T_{i+1}^{-}) \bLa_{q_{i}}\bx_{o}(T_{i+1}^{-}) + \beta^2 \bx_{c}^T(T_{i+1}^{-}) \bLa_{q_{i}}^{-1}\bx_{c}(T_{i+1}^{-}) \\
		&= V(\bx(T_{i+1}^{-}),\hat{\bx}(T_{i+1}^{-})),
		\end{split}
		\]
		i.e.,  \eqref{jumpVineq} holds. 
		\end{IEEEproof}
		\begin{IEEEproof}[Proof of Theorem \ref{theo:bound:step1}]
		From Lemma \ref{error:proof:lemma1} 
		it follows that
		\begin{multline*}
		V(\bx(s),\hat{\bx}(s)) - V(\bx(T_{i}),\hat{\bx}(T_i)) = \int_{T_i}^{s} \frac{\partial V(\bx_o(t), \bx_c(t))}{\partial t}dt\\ 
		\leqslant  4 \beta^2 \int_{T_i}^{s} \Vert \bu(t) \Vert_2^2dt- \int_{T_i}^{s} \Vert \by(t)-\hat{\by}(t) \Vert_2^2dt,  
		\end{multline*}
		and hence
		\begin{multline*}
		V(\bx(T_{i+1}^{-}),\hat{\bx}(T_{i+1}^{-})) - V(\bx(T_{i}),\hat{\bx}(T_i))\\ \leqslant  4 \beta^2 \int_{T_i}^{T_{i+1}} \Vert \bu(t) \Vert_2^2dt- \int_{T_i}^{T_{i+1}} \Vert \by(t)-\hat{\by}(t) \Vert_2^2dt.  
		\end{multline*}
		By Lemma \ref{error:proof:lemma2},  $V(\bx(T_{i+1}),\hat{\bx}(T_{i+1})) \le  V(\bx(T_{i+1}^{-}),\hat{\bx}(T_{i+1}^{-}))$ and hence
		\begin{multline*}
		V(\bx(T_{i+1}),\hat{\bx}(T_{i+1})) - V(\bx(T_{i}),\hat{\bx}(T_i))\\ \leqslant  4 \beta^2 \int_{T_i}^{T_{i+1}} \Vert \bu(t) \Vert_2^2dt- \int_{T_i}^{T_{i+1}} \Vert \by(t)-\hat{\by}(t) \Vert_2^2dt. 
		\end{multline*}
		By summing up the inequalities above,
		\[
		\begin{split}
		V(\bx(T_k),\hat{\bx}(T_k)) & - V(\bx(T_0),\hat{\bx}(T_0))\\ 
		&= \sum_{i=0}^{k-1}  V(\bx(T_{i+1}),\hat{\bx}(T_{i+1})) - V(\bx(T_{i}),\hat{\bx}(T_i))\\ 
		&\le  \sum_{i=0}^{k-1} 4 \beta^2\int_{T_i}^{T_{i+1}} \Vert \bu(t) \Vert_2^2dt\\
		&\qquad- \int_{T_i}^{T_{i+1}} \Vert \by(t)-\hat{\by}(t) \Vert_2^2dt \\
		& =4\beta^2 \int_{T_0}^{T_{k}} \Vert \bu(t) \Vert_2^2dt- \int_{T_0}^{T_{k}} \Vert \by(t)-\hat{\by}(t) \Vert_2^2dt.  
		\end{split}
		\]
		Using that $T_0=0$, $\bx(0)=0$, $\hat{\bx}(0)=0$, and $V(0,0)=0$ and  $V(\bx(T_k),\hat{\bx}(T_k)) \ge 0$, it follows that
		\[ 
		\begin{split}
		& 0 \le 4 \beta^2 \int_{0}^{T_{k}} \Vert \bu(t) \Vert_2^2dt- \int_{T_0}^{T_{k}} \Vert \by(t)dt-\hat{\by}(t) \Vert_2^2dt  \Leftrightarrow \\
		& \int_{T_0}^{T_{k}} \Vert \by(t)dt-\hat{\by}(t) \Vert_2^2dt \le 4\beta^2 \int_{0}^{T_{k}} \Vert \bu(t) \Vert_2^2dt.
		\end{split}
		\] 
		Since $\lim_{k\to\infty} T_k = \infty$, the statement of the theorem follows
		\end{IEEEproof}	
		
		\section{Numerical examples}
                \label{sect:num}
                
        \vspace{2mm}
                
		In this section, we analyze the practical applicability of the proposed MOR procedure. We consider a low-order artificial example represented by a linear hybrid systems with four subsystems. 
		
\begin{figure}[H]
\begin{center}
\includegraphics[scale=0.6]{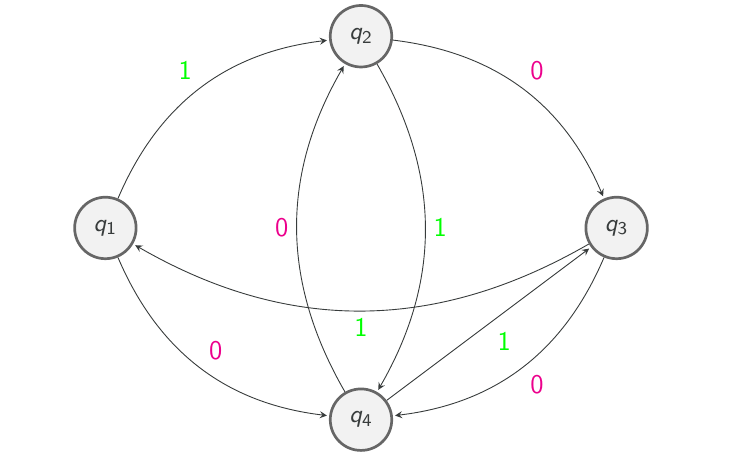}
\end{center}
\vspace{-4mm}
\caption{Directed graph representation of the state transition map.}
\label{fig:1}     
\end{figure}

		First, we characterize the discrete dynamics. The discrete state-transition map $\delta: \Omega \times \Gamma \rightarrow \Omega$ can be described in two ways, explicitly, i.e.:
		\begin{align*}
		\begin{cases}
		\textbf{Mode} \ \mathbf{q_1}: \ \ \delta(q_1,\textbf{\textcolor{magenta}{0}}) = q_4, \ \ \delta(q_1,\textbf{\textcolor{green}{1}}) = q_2, \\
		\textbf{Mode} \ \mathbf{q_2}: \ \ \delta(q_2,\textbf{\textcolor{magenta}{0}}) = q_3, \ \ \delta(q_2,\textbf{\textcolor{green}{1}}) = q_4, \\
		\textbf{Mode} \ \mathbf{q_3}: \ \ \delta(q_3,\textbf{\textcolor{magenta}{0}}) = q_4, \ \
		\delta(q_3,\textbf{\textcolor{green}{1}}) = q_1, \\
		\textbf{Mode} \ \mathbf{q_4}: \ \ \delta(q_4,\textbf{\textcolor{magenta}{0}}) = q_2, \ \ \delta(q_4,\textbf{\textcolor{green}{1}}) = q_3.
		\end{cases}
		\end{align*}
%
or using a directed graph, i.e. as in Fig.\;\ref{fig:1}.

\begin{figure}[H]
\begin{center}
\includegraphics[scale=0.27]{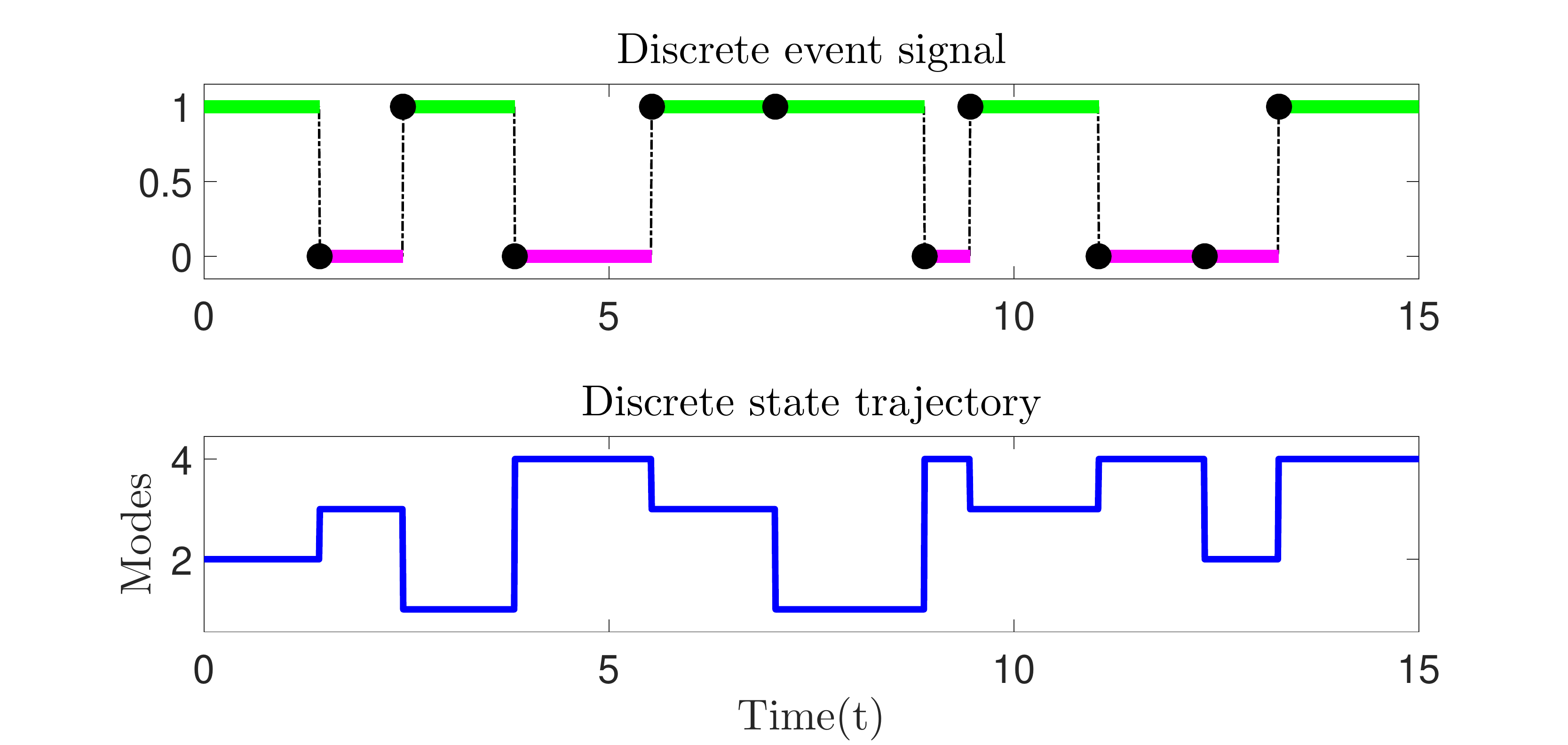}
\end{center}
\vspace{-4mm}
\caption{The discrete event signal $\gamma(t)$ (up) and the discrete state trajectory $\bq(t)$ (down).}
\label{fig:2}     
\vspace{-2mm}
\end{figure}

Next, we explicitly introduce the chosen discrete event signal $\gamma: \mathbb{R}_+ \rightarrow \Gamma$ and also the discrete state trajectory $\bq: \mathbb{R}_+ \rightarrow \Omega$
		\begin{align}\label{num_ex}
		\gamma(t) = \begin{cases} 1, \ t \in [0,T_1), \\ 0, \ t \in [T_1,T_2), \\ 1, \ t \in [T_2,T_3), \\ \ldots \\ 1, \ t \in [T_{10},T_{11}). \end{cases} \hspace{-2mm} \bq(t) = \begin{cases} q_2, \ t \in [0,T_1), \\ q_3, \ t \in [T_1,T_2), \\ q_1, \ t \in [T_2,T_3), \\ \ldots \\ q_4, \ t \in [T_{10},T_{11}), \end{cases} \hspace{-2mm}
		\end{align}
with given $T_1,\dots,T_{11}$ (see Fig.\;\ref{fig:2}). Additionally, in Fig.\;\ref{fig:2}, we depict the two signals introduced in (\ref{num_ex}), i.e. $\gamma(t)$ and $\bq(t)$ as a function of time (the time interval for this application was chosen to be $[0,15]$ seconds).
		

Finally, we proceed to the description of the continuous dynamics. Hence,	the system matrices $(A_q,B_q,C_q), 1 \leq q \leq 3$ corresponding to the linear hybrid system under consideration are written as follows:		
		
		\small
		\begin{align*}
		A_1 &= \left[\begin{array}{ccc} -1 & 0 & 0\\ 0 & -3 & 0\\ 0 & 0 & -4 \end{array}\right], \  A_2 = \left[\begin{array}{cc} -2 & 0\\ 0 & -1 \end{array}\right], \\ A_3 &= \left[\begin{array}{ccc} -3 & 0 & 0\\ 0 & -1 & 0\\ 0 & 0 & -2 \end{array}\right]
		, \ A_4 = \left[\begin{array}{cc} 1 & 0\\ 0 & \frac{1}{2} \end{array}\right], \\
		B_1 &= \left[\begin{array}{c} 1\\ -1\\ 1 \end{array}\right],  \ B_2 = \left[\begin{array}{c} 1\\ 1 \end{array}\right], \ \ B_3 = \left[\begin{array}{c} 1\\ 1\\ 3 \end{array}\right]
		\ , \\ B_4 &=  \left[\begin{array}{c} 2 \\ -1 \end{array}\right], \
		C_1 = \left[\begin{array}{ccc} 1 & -1 & 1 \end{array}\right], \ \ C_2 = \left[\begin{array}{cc} 1 & \frac{3}{2} \end{array}\right], \\ C_3 &= \left[\begin{array}{ccc} 1 & 1 & 1 \end{array}\right]
		, \ \ C_4 = \left[\begin{array}{cc} 2 & 1 \end{array}\right].
	\end{align*}

Additionally, the reset maps are given by the following matrices
		
		\begin{align*}
		M_{4,0,1} &= \frac{1}{\tau} \left[\begin{array}{ccc} 0 & 0 & -1\\ 0 & \frac{1}{2} & 0 \end{array}\right], \ \ M_{2,1,1}=   \frac{1}{\tau} \left[\begin{array}{ccc} 0 & 1 & 0\\ 1 & 0 & 0 \end{array}\right], \\ M_{3,0,2} &=  \frac{1}{\tau} \left[\begin{array}{cc} 0 & 1\\ 1 & 0\\ 0 & 0 \end{array}\right]
		,  \ \ M_{4,1,2} =  \frac{1}{\tau} \left[\begin{array}{cc} -1 & 1\\ 0 & 1 \end{array}\right], \\
		M_{4,0,3} &=   \frac{1}{\tau} \left[\begin{array}{ccc} 0 & 0 & 1\\ 0 & 0 & 0  \end{array}\right], \ \  M_{1,1,3} =  \frac{1}{\tau} \left[\begin{array}{ccc} 1 & -1 & 0\\ 0 & 0 & 1\\ 0 & -1 & 0 \end{array}\right], \\  M_{2,0,4} &=  \frac{1}{\tau} \left[\begin{array}{cc} -1 & 0\\ 0 & - \frac{1}{2} \end{array}\right], \ \   M_{3,1,4} =  \frac{1}{\tau} \left[\begin{array}{cc} -1 & 0\\ 1 & 0\\ 0 & \frac{1}{2} \end{array}\right].
		\end{align*}
		\normalsize
	
In the definition of the reset maps, one can observe that the scale $\tau >0$ is used. More precisely, in what follows, the value $\tau = 3$ was chosen for performing the numerical computations.
	
\begin{figure}[H]
		\begin{center}
		\includegraphics[scale=0.27]{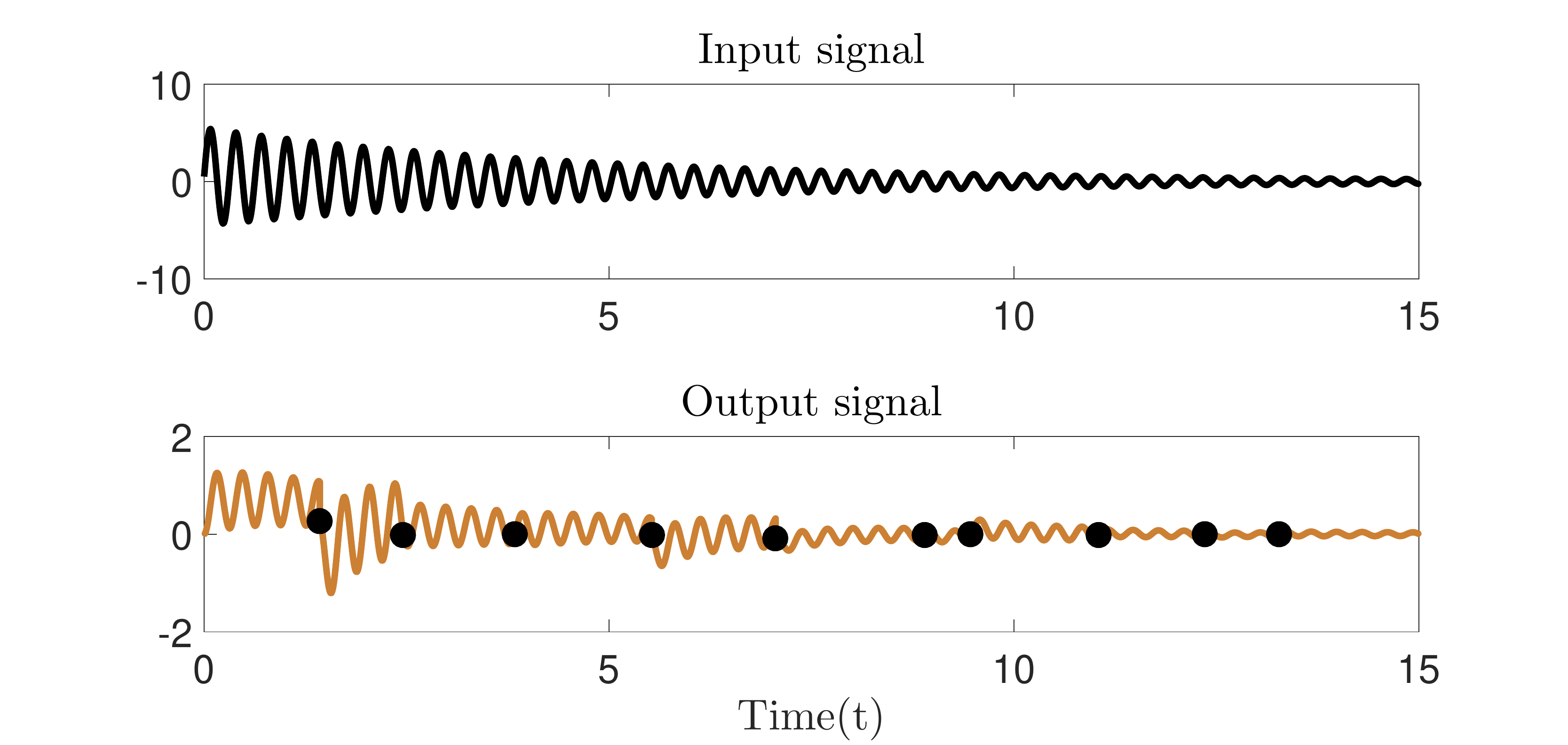}
		\end{center}
		\vspace{-4mm}
		\caption{The control input $\bu(t)$ (up) and the observed output $\by(t)$ (down).}
		\label{fig:3}     
		\vspace{-2mm}
		\end{figure}

We perform a time-domain simulation by using as continuous control input, the function $\bu(t) = 5 \sin(20t) e^{-t/5}+0.5 e^{-t/2}$. In Fig.\;\ref{fig:3}, we depict both the control input $\bu(t)$ and the observed output $\by(t)$ (as introduced in (\ref{petreczky:lti:hyb}))
		
		The next step is to find appropriate Gramians to be used in the balanced truncation procedure. We start by first computing the observability  Gramians.

		We are looking for positive definite matrices that satisfy the conditions in (\ref{obs:gram:eq}). Hence, for each mode, we explicitly state the corresponding LMIs:
		
		\vspace{4mm}

		\begin{itemize}
		
		\item Mode 1: $\begin{cases} A_1^T\cQ_1+\cQ_1 A_1+C_1^TC_1 < 0, \\
			M_{4,0,1}^T\cQ_{4} M_{4,0,1} - \cQ_{1} \leqslant 0,  \\
			M_{2,1,1}^T\cQ_{2} M_{2,1,1} - \cQ_{1} \leqslant 0 .
		\end{cases}$
		
		\item Mode 2: $\begin{cases} A_2^T\cQ_2+\cQ_2 A_2+C_2^TC_2 < 0, \\
			M_{3,0,2}^T\cQ_{3} M_{3,0,2} - \cQ_{2} \leqslant 0,  \\
			M_{4,1,2}^T\cQ_{4} M_{4,1,2} - \cQ_{2} \leqslant 0 .
		\end{cases}$

		\item Mode 3: $\begin{cases} A_3^T\cQ_3+\cQ_3 A_3+C_3^TC_3 < 0, \\
			M_{4,0,3}^T\cQ_{4} M_{4,0,3} - \cQ_{3} \leqslant 0,  \\
			M_{1,1,3}^T\cQ_{1} M_{1,1,3} - \cQ_{3} \leqslant 0 .
		\end{cases}$
		
		\item Mode 4: $\begin{cases} A_4^T\cQ_4+\cQ_4 A_4+C_4^TC_4 < 0, \\
			M_{2,0,4}^T\cQ_{2} M_{2,0,4} - \cQ_{4} \leqslant 0,  \\
			M_{3,1,4}^T\cQ_{3} M_{3,1,4} - \cQ_{4} \leqslant 0 .
		\end{cases}$
		
		\end{itemize}
		
		\vspace{4mm}
		
		It is to be remarked that, for $\tau = 1$, the above systems of LMIs could not be solved (by means of the optimization software provided in \cite{yalmip} and \cite{sedumi}). Nevertheless, when choosing $\tau = 3$, we were able to find a valid solution, i.e. a collection of positive definite matrices $\{\cQ_1, \cQ_2, \cQ_3, \cQ_4 \}$. More precisely, we could find:
		\begin{align*}
		\cQ_1 &= \left[ \begin{matrix}
		3.2662 &  -0.1118  &  0.0733 \\
		-0.1118 &   1.7564 &  -0.0693 \\
		0.0733 &  -0.0693 &   1.4755
		\end{matrix} \right], \\ \cQ_2 &=   \left[ \begin{matrix}      2.4546  & -0.0023 \\
		-0.0023  &  4.0827 \end{matrix} \right], \\
		\cQ_3 &=  \left[ \begin{matrix}
		1.7873  & -0.0041  &  0.0752 \\
		-0.0041  &  3.4766  &  0.1468 \\
		0.0752  &  0.1468  &  2.4182
		\end{matrix} \right], \\
		\cQ_4 &=  \left[ \begin{matrix}     3.9745  &  0.6789 \\
		0.6789  &  4.6925 
		\end{matrix} \right].
		\end{align*}
		
		Next, we need to find positive definite matrices $\cP_i$ that satisfy the conditions in (\ref{contr:gram:eq}). For each mode, we will state the corresponding LMIs:
		\vspace{2mm}
		\begin{itemize}
		\item Mode 1: $\begin{cases} A_1 \cP_1+\cP_1 A_1^T+B_1B_1^T < 0, \\
			M_{1,1,3} \cP_{3} M_{1,1,3}^T - \cP_{1} \leqslant 0,  \\
		\end{cases}$
		
		\item Mode 2: $\begin{cases} A_2 \cP_2+\cP_2 A_2^T+B_2 B_2^T < 0, \\
			M_{2,0,4} \cP_{4} M_{2,0,4}^T - \cP_{2} \leqslant 0, \\
			M_{2,1,1} \cP_{1} M_{2,1,1}^T - \cP_{2} \leqslant 0.
		\end{cases}$
		
		\item Mode 3: $\begin{cases} A_3 \cP_3+\cP_3 A_3^T+B_3 B_3^T < 0, \\
			M_{3,0,2} \cP_{2} M_{3,0,2}^T - \cP_{3} \leqslant 0,  \\
			M_{3,1,4} \cP_{4} M_{3,1,4}^T - \cP_{3} \leqslant 0 .
		\end{cases}$
		
		\item Mode 4: $\begin{cases} A_4 \cP_4+\cP_4 A_4^T+B_4 B_4^T < 0, \\
			M_{4,0,1} \cP_{1} M_{4,0,1}^T - \cP_{4} \leqslant 0,  \\
			M_{4,0,3} \cP_{3} M_{4,0,3}^T - \cP_{4} \leqslant 0,  \\
			M_{4,1,2} \cP_{2} M_{4,1,2}^T - \cP_{4} \leqslant 0 .
		\end{cases}$
				\vspace{2mm}
		\end{itemize}
	
		Again, for $\tau = 3$, we could find the following matrices
		\begin{align*}
		\cP_1 &= \left[ \begin{matrix}
		5.3173 &  -0.1332 &   0.3859 \\
		-0.1332 &   2.3055 &  -0.0914 \\
		0.3859 &  -0.0914  &  1.9288 
		\end{matrix} \right], \\ \cP_2 &=   \left[ \begin{matrix}  
		3.8471  &  0.1453 \\
		0.1453 &   5.3503
		\end{matrix} \right], \\
		\cP_3 &=  \left[ \begin{matrix}
		3.1234 &  -0.0344  &  0.3250 \\
		-0.0344 &   5.2759 &   0.5661 \\
		0.3250 &   0.5661  &  4.5523 \\
		\end{matrix} \right], \\ \cP_4 &=  \left[ \begin{matrix}  
		6.2062 &  -0.3344 \\
		-0.3344  &  7.4608
		\end{matrix} \right].
		\end{align*}
		
		Next, we present the Gramians in balanced representation, i.e. the diagonal matrices $\bLa_q$ from step 2 of Procedure 1.	
\begin{align*}
\bLa_1 &= \text{diag}(4.1894,  2.0184,  1.6542), \\ \bLa_2 &= \text{diag}(4.6754, 3.0703), \\
\bLa_3 &=   \text{diag}(4.3741,3.2543, 2.3291), \\ \bLa_4 &= \text{diag}( 5.9718,  4.8538).
\end{align*}

\begin{figure}[H]
		\begin{center}
		\includegraphics[scale=0.28]{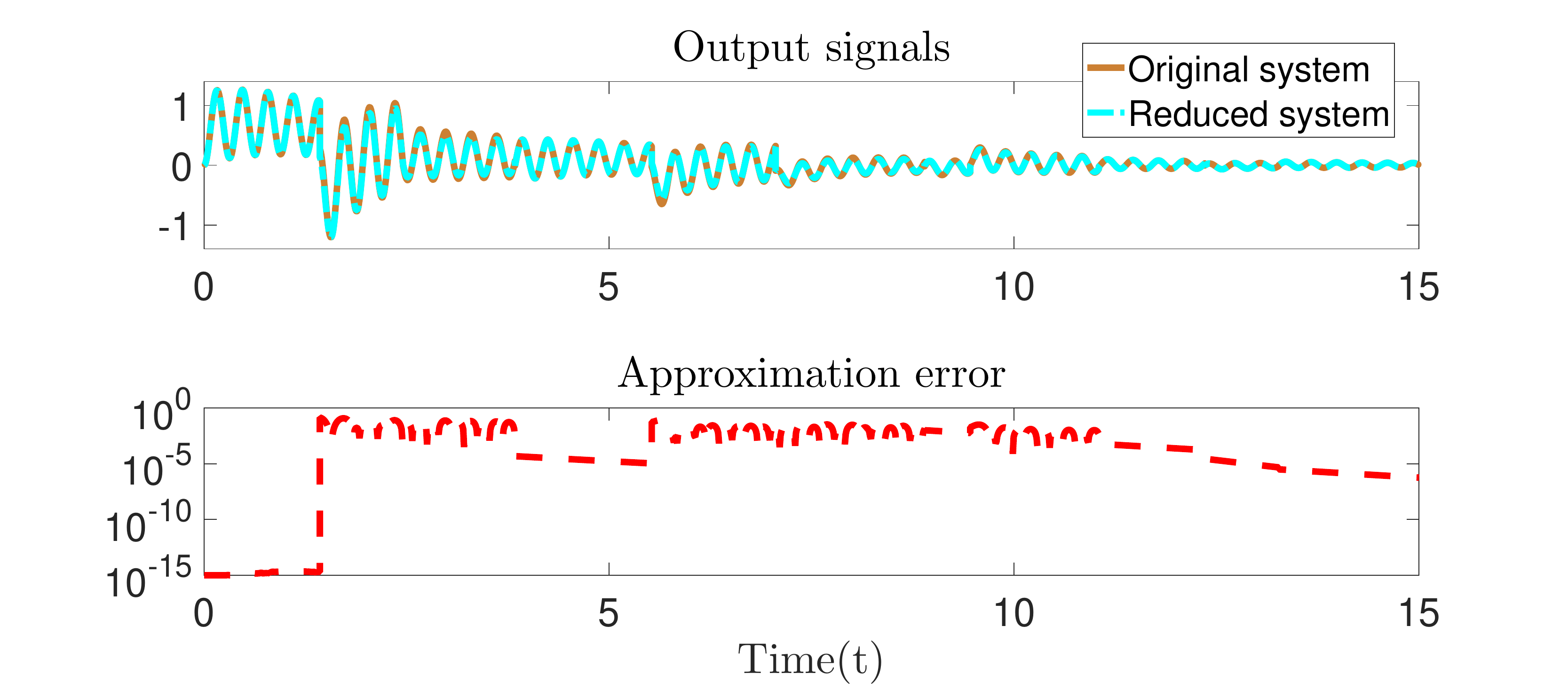}
		\end{center}
		\vspace{-4mm}
		\caption{The observed outputs for the original and reduced systems and the deviation between them (for the first choice of $r_k$'s).}
		\label{fig:4}     
		\vspace{-2mm}
		\end{figure}

		By choosing the reduction orders to be $r_1 = 2, r_2 = 2, r_3 = 2$ and $r_4=2$ (a dimension reduction is performed only for the first and third mode), we put together a reduced-order linear hybrid system. The time-domain simulation results are depicted in Fig.\;\ref{fig:4}.
	
	\vspace{-4mm}
		
%
%
		
\begin{figure}[H]
		\begin{center}
		\includegraphics[scale=0.28]{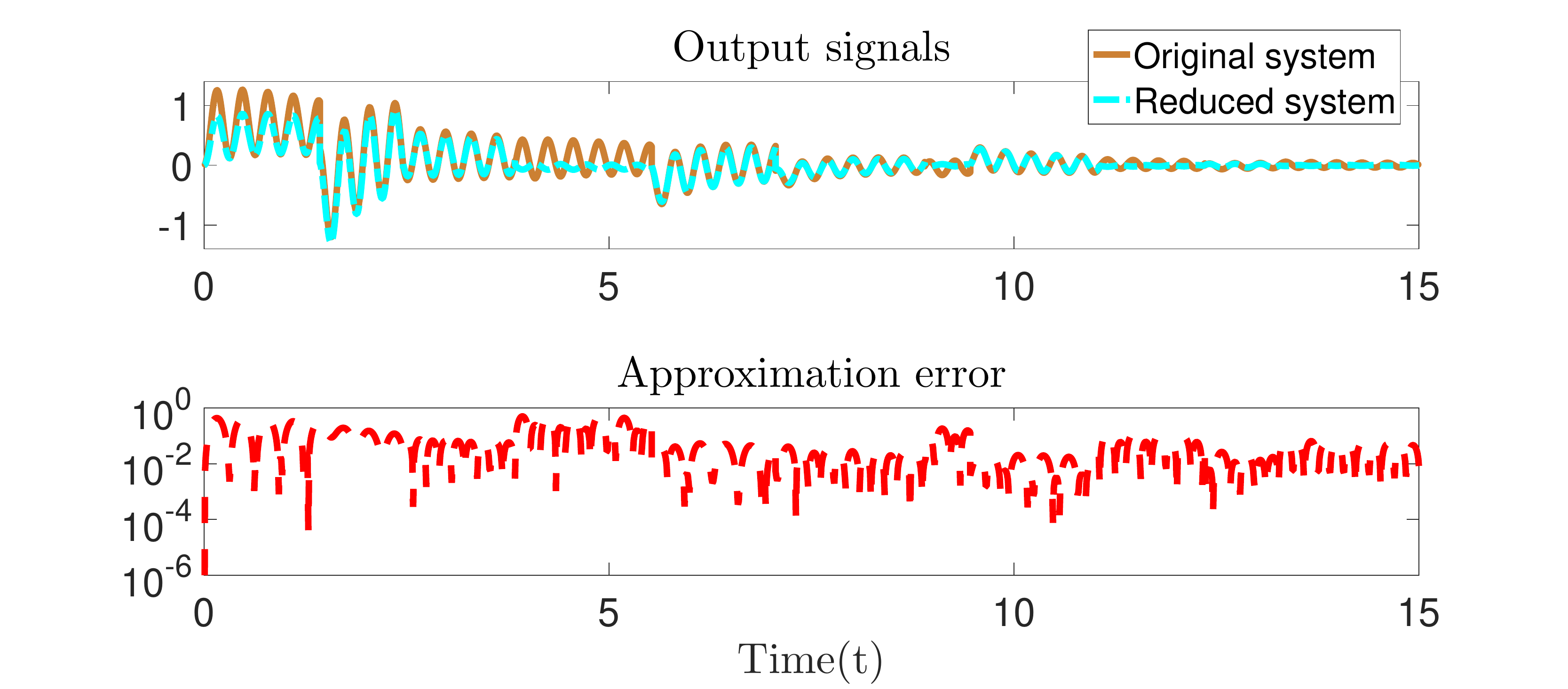}
		\end{center}
		\vspace{-4mm}
		\caption{The observed outputs for the original and reduced systems and the deviation between them (for the second choice of $r_k$'s).}
		\label{fig:5}     
		\vspace{-2mm}
		\end{figure}

		Next, we reduce the dimension of the systems corresponding to the second and forth modes as well. Hence, choose reduction orders  $r_1 = 2, r_2 = 1, r_3 = 2$ and $r_4=1$. The time-domain simulations results are depicted in  Fig.\;\ref{fig:5}.  
		
%
%

\section{Conclusion}

In this paper a balanced truncation procedure  for reducing linear hybrid systems was proposed. For each linear subsystem, specific Gramian matrices were computed by solving particular LMIs. An analytical error bound  in terms of singular values of the Gramians was also provided.   

We demonstrated the effectiveness of the procedure through a numerical example. Extensions that could be further developed include extending the proposed procedure to the case of hybrid systems with mild nonlinearities (such as systems with bilinear or stochastic behavior).

\appendix
\section{Technical proofs}
\label{tech:proofs}
\begin{IEEEproof}[Proof of Lemma \ref{lemma:stab}]
 Assume that $H$ is quadratically stable and assume that the positive definite matrices  $\{P_q\}_{q \in Q}$ satisfy \eqref{stab:gram:eq}. Then for suitable 
 $\gamma_q > 0$, $A_q^TP_q+P_qA_q < -\gamma_qI_{n_q}$. Note that $C_q^TC_q \le \mu_q I_{n_q}$ for a suitable $\mu_q > 0$. By taking 
 $\mu=\mathrm{min}\{ \frac{\gamma_q}{\mu_q}\}_{q  \in Q}$, it then follows that $A_q^TP_q+P_qA_q + \mu C_q^TC_q < 0$ from which it follows that
 $\cQ_q=\frac{1}{\mu} P_q$ is a generalized observability Gramian. 
 Similarly, by replacing $C_q^TC_q$ by $P_qB_qB_q^TP_q$ and repeating the argument above it follows that 
  $A_q^TP_q+P_qA_q + \mu P_qB_qB_q^TP_q < 0$ and by multiplying the latter LMI by $P_q^{-1}$ from right and left it follows that
  $A_qP_q^{-1}+P_q^{-1}A_q^T + \mu B_qB_q^T < 0$ from which, using the second equation of \eqref{stab:gram:eq} and \eqref{contr:gram:eq2:pf}
  it follows that $\cP_q=\frac{1}{\mu} P_q^{-1}$ is a generalized reachability Gramian. 
  Conversely, if $\{\cQ_q\}_{q \in Q}$ are generalized observability Gramians, then $A_q^T\cQ_q+\cQ_qA_q < -C_q^TC_q \le 0$ and hence $P_q=\cQ_q$ satisfy \eqref{stab:gram:eq}. 
  Similarly, if $\{\cP_q\}_{q \in Q}$ are generalized reachability Gramians, then by applying \eqref{contr:gram:eq2} with $u=0$ implies that $P_q=\cP_q^{-1}$, $q \in Q$ satisfy \eqref{stab:gram:eq}.
\end{IEEEproof}

\begin{IEEEproof}[Proof of Lemma \ref{observ:lemma}]
Let $\bx(t)$ be the corresponding solution to the \AHLS in (\ref{petreczky:ahlsd}), and  also introduce the function 
\begin{equation}\label{Vdef}
V(\bx(t)) = \begin{cases} \bx^T(t) \cQ_{q_0} \bx(t), \  t \in [0,t_1) \\ \bx^T(t) \cQ_{q_i} \bx(t),\  t \in [T_{i-1}, T_i), \ i \geqslant 2 \end{cases},
\end{equation}
where $T_i = \sum_{\ell=1}^{i} t_\ell$. By considering the uncontrolled case, the input function is considered to be $\bu(t) = 0, \ \forall t$. Using that $\frac{d \bx(t)}{dt} = A_{q_i} \bx(t)$, write the derivative of $V(t)$ from (\ref{Vdef}) for \textcolor{black}{$t \in [T_{i-1},T_i)$},
\begin{align*}
\frac{\partial V(\bx(t))}{\partial t} &= \frac{d \bx^T(t)}{d t} \cQ_{q_i} \bx(t) +  \bx^T(t) \cQ_{q_i}\frac{d \bx(t)}{d t} \\ &= \bx^T(t) \big{(} A_{q_i}^T \cQ_{q_i} + \cQ_{q_i} A_{q_i}  \big{)} \bx(t).
\end{align*}
By substituting the first inequality in (\ref{obs:gram:eq}) into the above relation, and using that $\by(t) = C_{q_i} \bx(t), \ t \in \textcolor{black}{[T_{i-1},T_i)}$, it follows that
\begin{equation}\label{Wderiv}
\frac{\partial V(\bx(t))}{\partial t} \leqslant - \by(t)^T \by(t).
\end{equation}
Introduce the following notation 
\begin{equation}\label{notationQ}
\bx(T_i^-) = \lim\limits_{t \nearrow T_i} \bx(t),   \ \ V(\bx(T_i^-)) = \lim\limits_{t \nearrow T_i} V(\bx(t)).
\end{equation}
By integrating the inequality (\ref{Wderiv}) from $T_{i-1}$ to \textcolor{black}{$t \in [T_{i-1},T_i)$}, it follows that
\begin{equation}\label{Wdiff}
V(\bx(t))-V(\bx(T_{i-1})) \leqslant -  \int\limits_{T_{i-1}}^t  \by(s)^T \by(s) ds.
\end{equation}
Using that $\bx(T_i) = M_{q_{i+1},\gamma,q_{i}} \bx(T_i^-)$, write
\begin{equation}\label{WTiplus}
V(\bx(T_i)) =  \bx^T(T_i^-) M_{q_{i+1},\gamma,q_{i}}^T \cQ_{q_{i+1}} M_{q_{i+1},\gamma,q_{i}} \bx(T_i^-).
\end{equation}
From the second inequality in (\ref{obs:gram:eq}), i.e. $M_{q_{i+1},\gamma,q_{i}}^T \cQ_{q_{i+1}} M_{q_{i+1},\gamma,q_{i}} \leqslant \cQ_{q_i}$, write
\begin{equation}\label{WTiplus2}
\begin{split}
V(\bx(T_i)) &= \bx^T(T_i^-) M_{q_{i+1},\gamma,q_{i}}^T \cQ_{q_{i+1}} M_{q_{i+1},\gamma,q_{i}} \bx(T_i^-) \\ &\leqslant \bx^T(T_i^-)  \cQ_{q_i}  \bx(T_i^-).
\end{split}
\end{equation}
Therefore, from (\ref{notationQ}), it follows that
\begin{equation}\label{WTiplus4}
V(\bx(T_i)) \leqslant V(\bx(T_i^-)).
\end{equation}
Putting together the inequalities in (\ref{Wdiff}) and (\ref{WTiplus4}), it follows that
\begin{equation}\label{WTiplus5}
V(\bx(T_i)) - V(\bx(T_{i-1})) \leqslant  -  \int\limits_{T_{i-1}}^{T_i}  \by(s)^T \by(s) ds.
\end{equation}
Now using the convention $T_0 = 0$ and adding all the inequalities in (\ref{WTiplus5}), we obtain
\begin{align}\label{WTiplus6}
\sum\limits_{i=1}^\ell V(\bx(T_i)) - V(\bx(T_{i-1})) \leqslant  -  \sum\limits_{i=1}^\ell \int\limits_{T_{i-1}}^{T_i}  \by(s)^T \by(s) ds \nonumber \\ \Rightarrow V(\bx(T_\ell)) - V(\bx(0)) \leqslant - \int\limits_{0}^{T_\ell}  \by(s)^T \by(s) ds.
\end{align}
Since $V(\bx(T_\ell)) = \bx^T(T_\ell) \cQ_{q_{\ell+1}} \bx(T_\ell)  \geqslant 0$, from (\ref{WTiplus6}) it follows that, 
\begin{equation}\label{WTiplus7}
V(\bx(0))   \geqslant  \int\limits_{0}^{T_\ell}  \by(s)^T \by(s) ds, \ \ \forall \ \ell \geq 0.
\end{equation}
By using that $V(\bx(0)) = \bx(0)^T \cQ_{q_0} \bx(0)$, the result in Lemma \ref{observ:lemma} is hence proven.
\end{IEEEproof}

\begin{IEEEproof}[Proof of Lemma \ref{control:lemma}]
Recall that $\cP_q$ satisfies the first inequality in (\ref{contr:gram:eq}).
By multiplying this inequality with $\cP_q^{-1}$ both to the left and to the right, we write
\begin{equation}\label{ineqP2}
A_q^T \cP_q^{-1} +\cP_q^{-1} A_q+ \cP_q^{-1} B_q B_q^T  \cP_q^{-1} < \bfz.
\end{equation}
Let $\bx(t)$ be the corresponding solution to the \AHLS in (\ref{petreczky:ahlsd}), and  also introduce the function  
\color{black}
\begin{equation}\label{Vdef2}
W(\bx(t)) = \begin{cases} \bx^T(t) \cP_{q_1}^{-1} \bx(t), \  t \in [0,t_1), \\ \bx^T(t) \cP_{q_i}^{-1} \bx(t),\  t \in [T_{i-1}, T_i), \ i \geqslant 2 \end{cases}.
\end{equation}
\color{black}
Using that $\dot{\bx}(t) = A_{q_i} \bx(t) + B_{q_i} \bu(t)$ and the definition of $W(\bx(t))$ in (\ref{Vdef2}), for 
\textcolor{black}{$t \in [T_{i-1},T_i)$}, we have
\begin{align*}
&\frac{\partial W(\bx(t))}{\partial t}= \frac{d \bx^T(t)}{d t} \cP_{q_i} ^{-1} \bx(t) +  \bx^T(t) \cP_{q_i}^{-1} \frac{d \bx(t)}{d t}\\ &= \bx^T(t) \big{(} A_{q_i}^T \cP_{q_i}^{-1} + \cP_{q_i}^{-1} A_{q_i}  \big{)} \bx(t) + 2 \bx(t)^T \cP_{q_i}^{-1} B_{q_i} \bu(t),
\end{align*}
and by using the inequality in (\ref{ineqP2}), it follows that
\begin{align}
\frac{\partial W(\bx(t))}{\partial t} &\leqslant -\bx(t)^T \cP_{q_i}^{-1} B_{q_i} B_{q_i}^T \cP_{q_i}^{-1} \bx(t) \nonumber + 2 \bx(t)^T \cP_{q_i}^{-1} B_{q_i} \bu(t) \nonumber \\ &= -\Vert B_{q_i}^T \cP_{q_i}^{-1} \bx(t) - \bu(t) \Vert_2^2 + \bu(t)^T \bu(t). 
\end{align}
Hence, the following inequality holds as,
\color{black}
\begin{equation}\label{Vineq}
\frac{\partial W(\bx(t))}{\partial t} \leqslant \bu(t)^T \bu(t), \ \ t \in [T_{i-1},T_i).
\end{equation}
\color{black}
Using (\ref{Vineq}) and integrating from $T_{i-1}$ to t, we obtain
\begin{equation}\label{Wdiff2}
W(\bx(t))-W(\bx(T_{i-1})) \leqslant  \int\limits_{T_{i-1}}^t \bu^T(s) \bu(s) ds .
\end{equation}
Using that $\bx(T_i) = M_{q_{i+1},\gamma,q_{i}} \bx(T_i^-)$, write
\begin{equation}\label{WTiplus11}
W(\bx(T_i)) =  \bx^T(T_i^-) M_{q_{i+1},\gamma,q_{i}}^T \cP_{q_{i+1}}^{-1} M_{q_{i+1},\gamma,q_{i}} \bx(T_i^-).
\end{equation}
From the second inequality in (\ref{contr:gram:eq2}), one can directly derive that $M_{q_{i+1},\gamma,q_{i}}^T \cP_{q_{i+1}}^{-1} M_{q_{i+1},\gamma,q_{i}} \leqslant \cP_{q_i}^{-1}$. Then, 
\begin{equation}\label{WTiplus22}
\begin{split}
W(\bx(T_i)) &= \bx^T(T_i^-) M_{q_{i+1},\gamma,q_{i}}^T \cP_{q_{i+1}}^{-1} M_{q_{i+1},\gamma,q_{i}} \bx(T_i^-) \\ &\leqslant \bx^T(T_i^-)  \cP_{q_i}^{-1}  \bx(T_i^-).
\end{split}
\end{equation}
Therefore, it follows that $W(\bx(T_i)) \leqslant W(\bx(T_i^-))$, where $W(\bx(T_i^-)) = \lim\limits_{t \uparrow T_i} W(\bx(t))$ for $i>0$ and $W(\bx(0^-)) = W(\bx(0))$.

By combining this inequality with the inequality in (\ref{Wdiff2}), one can write
\begin{align}\label{Wineq_cont1}
& W(\bx(T_i^-)) - W(\bx(T_{i-1}^-)) \leqslant  \int\limits_{T_{i-1}}^{T_i} \bu^T(s) \bu(s) ds \nonumber~~\Rightarrow \\ &  \sum\limits_{i=1}^\ell W(\bx(T_i^-)) - W(\bx(T_{i-1}^-)) \leqslant    \sum\limits_{i=1}^\ell  \int\limits_{T_{i-1}}^{T_i}  \bu^T(s) \bu(s) ds ~~\Rightarrow \nonumber \\ &   W(\bx(T_\ell^-)) - W(\bx(0^-)) \leqslant  \int\limits_{0}^{T_\ell}  \bu^T(s) \bu(s) ds.
\end{align}
Since $\bx(0) = \bfz$, it follows that $W(\bx(0^-)) = \bfz$. Also, from the definition of the function W, it is clear that $W(\bx(T_\ell^-)) = \bx^T(T_\ell^-) \cP_{q_\ell}^{-1} \bx(T_\ell^-)$. Hence, from (\ref{Wineq_cont1}), we directly conclude that
\begin{equation}\label{ineqP_Tl}
\bx^T(T_\ell^-) \cP_{q_\ell}^{-1} \bx(T_\ell^-) \leqslant \int\limits_{0}^{T_\ell}  \bu^T(s) \bu(s) ds, \ \forall \ell \geqslant 1,
\end{equation}
which proves the result in Lemma \ref{control:lemma}. 
\end{IEEEproof}

\begin{IEEEproof}[Proof of Lemma \ref{balanced:lemma1}]
It is easy to see that $\bS_q^{T} \bLa_q\bS_q=\cQ_q$ and $\bS_q^{-1}\bLa_q\bS_q^{-T}=\cP_q$. 
From  $\bS_q^{T}\bLa_q\bS_q=\cQ_q$ it follows that 
\[
\begin{split}
& \bar{A}_q^T\bLa_q+\bLa_q\bar{A}_q+\bar{C}_q^T\bar{C}_q < 0,  \\
& \forall \gamma \in \Gamma, q^{+}=\delta(q,\gamma): 
\bar{M}_{q^+,\gamma,q}^T\bLa_{q^+} \bar{M}_{q^+,\gamma,q} - \bLa_{q} \leqslant 0,
\end{split}
\]
which means that $\{\bLa_q\}_{q \in Q}$ are generalized observability Gramians of $\bar{H}$.
Indeed, by using \eqref{eq:bal1},
\[
\begin{split}
& \bar{A}_q^T\bLa_q+\bLa_q\bar{A}_q+\bar{C}_q^T\bar{C}_q=
\bS_q^{-T} A_q^T \bS_q^T \bLa_q + \bLa_q \bS_q A_q \bS_q^{-1}\\ &+ \bS_q^{-T}C_q^TC_q \bS_q^{-1}  
= \bS_q^{-T}(A_q^T \underbrace{\bS_q^T\bLa_q \bS_q}_{=\cQ_q} + \underbrace{\bS_q^{T}\bLa_q\bS_q}_{=\cQ_q} A_q + C_q^TC_q)\bS_q^{-1} \\
& = \bS_q^{-T}(A_q^T \cQ_q + \cQ_q A_q + C_q^TC_q)\bS_q^{-1}.
\end{split}
\]
Since $A_q^T \cQ_q + \cQ_q A_q + C_q^TC_q < 0$, it follows that
\begin{align*}
& \bar{A}_q^T\bLa_q+\bLa_q\bar{A}_q+\bar{C}_q^T\bar{C}_q \\ &=
\bS_q^{-T}(A_q^T \cQ_q + \cQ_q A_q + C_q^TC_q)\bS_q^{-1} < 0. 
\end{align*}
Similarly,
\[
\begin{split}
& \bar{M}_{q^+,\gamma,q}^T\bLa_{q^+} \bar{M}_{q^+,\gamma,q} - \bLa_{q} = \bS_{q}^{-T}M_{q^+,\gamma,q}^T\bS_{q^+}^T\bLa_{q^+} \bS_{q^{+}}M_{q^+,\gamma,q}\bS_q^{-1} \\ & - \bLa_{q}
=\bS_q^{-T}(M_{q^+,\gamma,q}^T\underbrace{\bS_{q^+}^T\bLa_{q^+} \bS_{q^{+}}}_{=\cQ_{q^{+}}}M_{q^+,\gamma,q} - \underbrace{\bS_q^T\bLa_{q}\bS_q}_{=\cQ_q})\bS_q^{-1}\\
& = \bS_q^{-T}(M_{q^+,\gamma,q}^T\cQ_q^{+}M_{q^+,\gamma,q} - \cQ_q)\bS_q^{-1}.
\end{split}
\]
Since $(M_{q^+,\gamma,q}^T\cQ_{q^{+}}M_{q^+,\gamma,q} - \cQ_q) \leq 0$, it then follows that
\begin{align*}
&\bar{M}_{q^+,\gamma,q}^T\bLa_{q^+} \bar{M}_{q^+,\gamma,q} - \bLa_{q}\\ & = \bS_q^{-T}(M_{q^+,\gamma,q}^T\cQ_q^{+}M_{q^+,\gamma,q} - \cQ_q)\bS_q^{-1}	\leq 0.
\end{align*}
The proof that $\{\bLa_q\}_{q \in Q}$ are generalized reachability Gramians is similar to the proof above.
\end{IEEEproof}

\begin{IEEEproof}[Proof of Lemma \ref{balanced:lemma}]
We will show that $\{\hat{\bLa}_q\}_{q \in Q}$ are observability Gramians, the proof that it is a reachability Gramian is completely analogous.
The claim of the lemma on quadratic stability of $\hat{H}$ follows from Lemma \ref{lemma:stab}.
First, we show that $\hat{A}_q^T\hat{\bLa}_q+\hat{\bLa}_q\hat{A}_q+\hat{C}_q^T\hat{C}_q < 0$ for all $q \in Q$.
If $r_q=n_q$, then $(\bar{A}_q,\bar{B}_q,\bar{C}_q,\bLa_q)=(\hat{A}_q,\hat{B}_q,\hat{C}_q,\hat{\bLa}_q)$, and as by  Lemma \ref{balanced:lemma1} it follows that $\{\bLa_q\}_{q \in Q}$ is a o observability Gramian, $\hat{A}_q^T\hat{\bLa}_q+\hat{\bLa}_q\hat{A}_q+\hat{C}_q^T\hat{C}_q < 0$  holds.
If $r_q < n_q$, then   
\small
\begin{equation}
\label{balanced:lemma:pf:eq1}
\begin{split}
&\bar{A}_q^T\bLa_q+\bLa_q\bar{A}_q+\bar{C}_q^T\bar{C}_q=
\begin{bmatrix} (\bar{A}_q^{11})^T & (\bar{A}_q^{21})^T \\ (\bar{A}_q^{12})^T & (\bar{A}_q^{22})^T \end{bmatrix} \begin{bmatrix} \hat{\bLa}_q & 0 \\ 0 & \beta_q \end{bmatrix}\\
&+
\begin{bmatrix} \hat{\bLa}_q & 0 \\ 0 & \beta_q \end{bmatrix}
\begin{bmatrix} \bar{A}_q^{11} & \bar{A}_q^{12} \\ \bar{A}_q^{21} & \bar{A}_q^{22} \end{bmatrix}  + 
\begin{bmatrix} (\bar{C}_q^{1})^T \\  (\bar{C}_q^{2})^T \end{bmatrix}\begin{bmatrix} \bar{C}_q^{1} & \bar{C}_q^{2} \end{bmatrix} = \\
&\begin{bmatrix} 
(\bar{A}_q^{11})^T\hat{\bLa}_q+ \hat{\bLa}_q\bar{A}_q^{11} + (\bar{C}_q^{1})^T\bar{C}_q^{1} & \star \\ \star & \star \end{bmatrix}=\begin{bmatrix} \hat{A}_q^T\hat{\bLa}_q+\hat{\bLa}_q\hat{A}_q+\hat{C}_q^T\hat{C}_q & \star \\ \star & \star \end{bmatrix}.
\end{split}
\end{equation}
\normalsize
From Lemma \ref{balanced:lemma1} it follows that $\{\bLa_q\}_{q \in Q}$ are observability Gramians, and thus
$\bar{A}_q^T\bLa_q+\bLa_q\bar{A}_q+\bar{C}_q^T\bar{C}_q < 0$ holds. This implies that the left-upper $r_q \times r_q$ block of $\bar{A}_q^T\bLa_q+\bLa_q\bar{A}_q+\bar{C}_q^T\bar{C}_q$,
which equals $\hat{A}_q^T\hat{\bLa}_q+\hat{\bLa}_q\hat{A}_q+\hat{C}_q^T\hat{C}_q$ is also negative definite.

Next, we show that
\begin{equation}
\label{balanced:lemma:pf:eq2}
\hat{M}_{q^+,\gamma,q}^T\hat{\bLa}_{q^+} \hat{M}_{q^+,\gamma,q} - \hat{\bLa}_{q} \leqslant 0. 
\end{equation}
If $r_q=n_q, r_{q^{+}}=n_{q^{+}}$, then $\hat{M}_{q^+,\gamma,q}=\bar{M}_{q^{+},\gamma,q}$, 
$\bLa_q=\hat{\bLa}_q$, $\bLa_{q^+}=\hat{\bLa}_q^{+}$, and as $\bar{M}_{q^+,\gamma,q}^T\bLa_{q^+} \bar{M}_{q^+,\gamma,q} - \bLa_{q} \leq 0.$
\eqref{balanced:lemma:pf:eq2} follows. 
For the other cases, we proceed to prove that  
\[
\hat{M}_{q^+,\gamma,q}^T\hat{\bLa}_{q^+} \hat{M}_{q^+,\gamma,q} - \hat{\bLa}_{q} = \begin{bmatrix} D & \star \\ \star & \star \end{bmatrix},
\]
where the matrix $D$ is such that
\[ D \ge \hat{M}_{q^+,\gamma,q}^T\hat{\bLa}_{q^+} \hat{M}_{q^+,\gamma,q} - \hat{\bLa}_{q}. \]
If this is the case, then from \eqref{balanced:lemma:pf:eq2} it follows that $D \leq 0$, from which it follows that
$\hat{M}_{q^+,\gamma,q}^T\hat{\bLa}_{q^+} \hat{M}_{q^+,\gamma,q} - \hat{\bLa}_{q} \le 0$.
Consider the case when $r_{q^{+}} < n_{q^{+}}$ and $r_q < n_q$. 
\small
\begin{align*}
\begin{split}
&\bar{M}_{q^+,\gamma,q}^T\bLa_{q^+} \bar{M}_{q^+,\gamma,q} - \bLa_{q} = 
\begin{bmatrix} (\bar{M}_{q^+,\gamma,q}^{11})^T & (\bar{M}_{q^+,\gamma,q}^{21})^T \\ (\bar{M}_{q^+,\gamma,q}^{12})^T & (\bar{M}_{q^{+},\gamma,q}^{22})^T
\end{bmatrix} \\
&\begin{bmatrix} \hat{\bLa}_{q^+} & 0 \\ 0 & \beta_{q^+} \end{bmatrix}
\begin{bmatrix} \bar{M}_{q^+,\gamma,q}^{11} & \bar{M}_{q^+,\gamma,q}^{12} \\ \bar{M}_{q^+,\gamma,q}^{21} & \bar{M}_{q^+,\gamma,q}^{22}
\end{bmatrix}
- \begin{bmatrix} \hat{\bLa}_q & 0 \\ 0 & \beta_{q} \end{bmatrix}= \\
&\begin{bmatrix} 
\underbrace{(\bar{M}_{q^+,\gamma,q}^{11})^T\hat{\bLa}_{q^+} \bar{M}_{q^+,\gamma,q}^{11}+\beta_{q^{+}}(\bar{M}_{q^+,\gamma,q}^{21})^T\bar{M}_{q^+,\gamma,q}^{21} - \hat{\bLa}_q}_{=D} & \star \\ \star & \star
\end{bmatrix}. 
\end{split}
\end{align*}
\normalsize
In this case, since $\beta_{q^{+}}(\bar{M}_{q^+,\gamma,q}^{21})^T\bar{M}_{q^+,\gamma,q}^{21} \ge 0$, it follows that
\[ 
\begin{split}
&D=(\bar{M}_{q^+,\gamma,q}^{11})^T\hat{\bLa}_{q^+}\bar{M}_{q^+,\gamma,q}^{11}- \hat{\bLa}_q+\beta_{q^{+}}(\bar{M}_{q^+,\gamma,q}^{21})^T\bar{M}_{q^+,\gamma,q}^{21} \\ &\geq 
(\bar{M}_{q^+,\gamma,q}^{11})^T\hat{\bLa}_{q^{+}}\bar{M}_{q^+,\gamma,q}^{11}- \hat{\bLa}_q =
\hat{M}_{q^+,\gamma,q}^T\hat{\bLa_q^+}\hat{M}_{q^+,\gamma,q}- \hat{\bLa}_q.
\end{split}   
\]
If $r_{q^{+}}=n_{q^{+}}$ but $r_q < n_q$, then $\hat{\bLa}_{q^+}=\bLa_{q^{+}}$, and 
\small
\[
\begin{split}
&\bar{M}_{q^+,\gamma,q}^T\bLa_{q^+} \bar{M}_{q^+,\gamma,q} - \bLa_{q} = 
\begin{bmatrix} (\bar{M}_{q^+,\gamma,q}^{11})^T \\ (\bar{M}_{q^+,\gamma,q}^{12})^T \end{bmatrix}
\hat{\bLa}_{q^+} 
\begin{bmatrix} \bar{M}_{q^+,\gamma,q}^{11} & \bar{M}_{q^+,\gamma,q}^{12}
\end{bmatrix} \\
&- \begin{bmatrix} \hat{\bLa}_q & 0 \\ 0 & \beta_{q} \end{bmatrix} =
\begin{bmatrix} 
\underbrace{(\bar{M}_{q^+,\gamma,q}^{11})^T\hat{\bLa}_{q^+}\bar{M}_{q^+,\gamma,q}^{11} - \hat{\bLa}_q}_{=D} & \star \\ \star & \star
\end{bmatrix}. 
\end{split}
\]
\normalsize
In this case, $D=\hat{M}_{q^+,\gamma,q}^T\hat{\bLa}_{q^+} \hat{M}_{q^+,\gamma,q}- \hat{\bLa}_q$.
Finally, if
$r_{q^{+}} < n_{q^{+}}$ but $r_q = n_q$, then $\hat{\bLa}_{q}=\bLa_{q}$, and 
\small
\[
\begin{split}
&\bar{M}_{q^+,\gamma,q}^T\bLa_{q^+} \bar{M}_{q^+,\gamma,q} - \bLa_{q} = 
\begin{bmatrix} (\bar{M}_{q^+,\gamma,q}^{11})^T & (\bar{M}_{q^+,\gamma,q}^{21})^T\end{bmatrix} \\
&\begin{bmatrix} \hat{\bLa}_{q^+} & 0 \\ 0 & \beta_{q^{+}} \end{bmatrix} 
\begin{bmatrix} \bar{M}_{q^+,\gamma,q}^{11} \\ \bar{M}_{q^+,\gamma,q}^{21}
\end{bmatrix}
- \hat{\bLa}_q \\
&=
(\bar{M}_{q^+,\gamma,q}^{11})^T\hat{\bLa}_{q^+}\bar{M}_{q^+,\gamma,q}^{11} - \hat{\bLa}_q +  \beta_{q^{+}}(\bar{M}_{q^+,\gamma,q}^{21})^T\bar{M}_{q^+,\gamma,q}^{21}=D
\end{split}
\]
\normalsize
and in this case $D \geq \hat{M}_{q^+,\gamma,q}^T\hat{\bLa}_{q^+}\hat{M}_{q^+,\gamma,q}- \hat{\bLa}_q$ since $\beta_{q^{+}}(\bar{M}_{q^+,\gamma,q}^{21})^T\bar{M}_{q^+,\gamma,q}^{21} \ge 0$.
\end{IEEEproof}


%
%

\ifCLASSOPTIONcaptionsoff
  \newpage
\fi



%
%
%


\bibliographystyle{IEEEtran}
\bibliography{HybridBT2019}

\end{document}